\documentclass[11pt]{amsart}
\usepackage{fullpage,amsfonts,amsmath,amscd,amssymb}
\theoremstyle{plain}

\newtheorem*{lem}{Lemma}
\newtheorem*{prop}{Proposition}
\newtheorem*{thm}{Theorem}
\newtheorem*{cor}{Corollary}

\theoremstyle{remark}
\newtheorem*{rem}{Remark}
\newtheorem*{rems}{Remarks}
\newcommand{\x}[4]{\mbox{Ext}_{#1}^{#2}({#3},{#4})}
\newcommand{\hm}[3]{\mbox{Hom}_{#1}({#2},{#3})}
\newcommand{\ep}{\epsilon}
\newcommand{\g}{\mathfrak{g}}

\newcommand{\C}{\mathbb{C}}

\newcommand{\UE}{U_{\epsilon}(\mathfrak{g})}  
\newcommand{\UM}{U_{\epsilon}^{\leq 0}}  
\newcommand{\UP}{U_{\epsilon}^{\geq 0}}  
\newcommand{\BARUP}{\overline{U_{\epsilon}^{\geq 0}}}   
\newcommand{\EO}{\mathcal{O}_{\epsilon}[G]}   
\newcommand{\BARUN}{\overline{U_{\epsilon}^{>0}}}
\newcommand{\bi}{b_{\varpi_i}}
\newcommand{\ci}{c_{\varpi_i}}

\title{The Ramifications of the centres: quantised function algebras at roots of unity}
\author{Kenneth A. Brown and Iain Gordon}
\address{Department of Mathematics, University of Glasgow, Glasgow G12 8QW}
\email{kab@maths.gla.ac.uk}
\email{igordon@Mathematik.Uni-Bielefeld.DE}
\thanks{The research
 of the first author was partially supported by NATO Grant 960250. The second author undertook
 part of the research for this paper while holding a William
Seggie Brown Postdoctoral Fellowship at the University of
Edinburgh, and part while supported by TMR grant ERB
FMRX-CT97-0100 at the University of Bielefeld. Both authors thank
Toby Stafford for helpful conversations.}
\begin{document}
\begin{abstract}
This paper continues the
study of quantised function algebras $\EO$ of a semisimple group
$G$ at an $\ell$th root of unity $\epsilon$. These algebras were
introduced by De Concini and Lyubashenko in 1994, and studied
further by
De Concini and Procesi and by Gordon, amongst others. Our main
purpose here is to
increase understanding of the finite dimensional factor
algebras $\EO (g)$, for $g \in G$. We determine the representation
type and block structure of these factors, and (for many $g$)
describe them up to isomorphism. A series of parallel results is
obtained for the quantised Borel algebras $\UP$ and $\UM$.
\end{abstract}
\maketitle
\section{Introduction}
\subsection{}
\label{I1} The first substantial study of the quantised function
algebra $\EO$ of the simply-connected semisimple group $G$ at the
$\ell$th root of unity $\epsilon$ appeared in \cite{declyu1}. It
was shown there that, in close analogy with the case of a generic
parameter \cite{jos}, the representation theory of $\EO$ is
stratified by double Bruhat cells in $G$. More precisely, $\EO$
contains a central sub-Hopf algebra isomorphic to
$\mathcal{O}[G]$, over which $\EO$ is a projective module of
constant rank $\ell^{\mathrm{dim}G}$. (In fact as we show below in
Proposition 2.2, $\EO$ is a free $\mathcal{O}[G]$-module.)
Every irreducible $\EO$-module is annihilated by a maximal ideal
$\mathfrak{m}_g$ of $\mathcal{O}[G]$ (where $g \in G$), so that
the finite dimensional representation theory of $\EO$ can
effectively be reduced to the study of the bundle of
$\ell^{\mathrm{dim}G}$-dimensional algebras $\EO (g):=
\EO/\mathfrak{m}_g \EO$, for $g \in G$.
The noncommutativity of the generic algebra $\mathcal{O}_q [G]$
induces a Poisson structure on $\mathcal{O}[G]$ (as in
\cite[Section 11]{decpro3}), and this is preserved by the group
$T$ of winding automorphisms of $\EO$ afforded by the
one-dimensional representations of $\EO$. (Here, $T$ is the
maximal torus in $G$.) It follows (see \cite[Section 9]{declyu1})
that if $g$ and $g'$ are in the same $T$-orbit of symplectic
leaves in $G$, then $\EO (g) \cong \EO (g')$ \cite[9.3]{declyu1}.
The $T$-orbits of symplectic leaves have been determined
\cite{hl1}: they are the double Bruhat cells
$$ X_{w_1,w_2} = B^+ w_1 B^+ \cap B^- w_2  B^-, $$
where $w_1, w_2 \in W$, the Weyl group of $G$, and $B^{+}$ and
$B^-$ are fixed Borel subgroups of $G$.
\subsection{}
\label{I2} The representation theory of $\EO (g)$ was further
studied in \cite{decpro49}, where it was shown that the
irreducible $\EO (g)$-modules are permuted transitively by the
winding automorphisms arising from the one-dimensional $\EO
(1)$-modules, and the number and dimension of irreducible $\EO
(g)$-modules was calculated - see Theorem \ref{simplefun}(b)(ii).
The complexity of $\EO (g)$ was determined in \cite{gor4}, and
hence the representation type of $\EO (g)$ was found in many, but
not all, cases - see Theorem \ref{tamewild}(b).
\subsection{}
\label{I3} The main purpose of this paper is to continue the
analysis of the algebras $\EO (g)$, our principal results in this
direction being listed below. For $w \in W$ let
$\ell(w)$ (respectively $s(w)$) denote the minimal length of an
expression for $w$ as a product of simple (respectively arbitrary)
reflections in $W$. Let $w_1,w_2 \in W$, let $g \in
X_{w_1,w_2}$ and set $w = w_2^{-1}w_1$. Let $N$ be the number of positive roots of
$G$, let $r$ be the rank of $G$, let $\varpi_1, \ldots , \varpi_r$
be a set of fundamental weights, and let $\mathfrak{S}(w_1,w_2) =
\{ i
: 1 \leq i \leq r, w_0w_1, w_0w_2 \in \mathrm{Stab}_W (\varpi_i)
\}$. By a \textit{multiply-edged Cayley graph of} $F$ we mean the
graph got from the usual Cayley graph $C$ of a group $F$ with
respect to a distinguished set of generators $X$ (possibly
including $1_F$) by assigning a positive integer $m_x$ to each $x
\in X$ and replacing each edge of $C$ corresponding to $x$ by
$m_x$ edges in the same direction.
\begin{itemize}
\item (Theorem \ref{azmor}) For $g$ in the fully Azumaya locus
(that is, for those algebras whose irreducible modules have the
maximal possible dimension $\ell^N$), a complete description of
$\EO (g)$ as a direct sum of matrix rings over a truncated
polynomial algebra.
\item (Theorem \ref{reptypethm2}) Determination
of the representation type of $\EO (g)$ in all cases - $\EO (g)$
has finite type if $\ell (w_1) + \ell (w_2) > 2N - 2$, and is wild
otherwise.
\item (Corollary \ref{oblocks}) Calculation of the number of
blocks of $\EO (g)$ and the structure of its quiver: the number of
blocks is $\ell^{\mathrm{card}\mathfrak{S}(w_1, w_2)}$ and the
quiver of each block is a multiply-edged Cayley graph of a certain elementary
abelian $\ell$-group of order $\ell^{r - s(w)
-\mathrm{card}\mathfrak{S}(w_1, w_2)}$.
\end{itemize}
\subsection{}
\label{I4} Here are some indications of the ingredients of the
proofs of the above results. This paper is a sequel to
\cite{brogor1}, in which a key result (reproduced here as Theorem
\ref{maxfactor}) shows in the present setting that if $g$ is a
fully Azumaya point of $G$ then $\EO/\mathfrak{m}_g \EO$ is a
complete matrix ring over $Z_g := Z(\EO)/\mathfrak{m}_g Z( \EO)$.
A particular case of the main result of \cite{brogoo} (Theorem
\ref{az} below) identifies the Azumaya points of $\EO$ with the
smooth points of its centre. Since $Z(\EO)$ is known thanks to
work of Enriquez \cite{enr}, $Z_g$ can be determined when $g$ lies
under a smooth point of $Z(\EO)$, so proving Theorem \ref{azmor}.
For the analysis of representation type, after the work of
\cite{gor4}, (at least for $\ell$ greater than the Coxeter number
$h$ of $G$), only the case (*) $\ell (w_1) + \ell (w_2) = 2N - 2$
remained to be dealt with. In this case (and without assuming
$\ell > h$) there are essentially three possibilities for the pair
$(w_1,w_2)$, two of which yield an Azumaya point which can thus be
disposed of thanks to Theorem \ref{azmor}. The third (where $w_1 =
w_2$ and $\ell (w_1) = N - 1$) is then analysed directly, and
shown always to involve the wild algebra
$\C[X,Y]/(X^{\ell},Y^{\ell})$. Finally, using deformation
arguments we remove the restriction to $\ell > h$ arising in
\cite{gor4} - the idea is that every algebra $\EO (g')$ with $g'
\in X_{w_1',w_2'}$ and $\ell (w_1') + \ell (w_2') < 2N - 2$ is a
degeneration of an algebra $\EO (g)$ for $g \in X_{w_1,w_2}$
satisfying (*). Then by a result of Gei{\ss} \cite{gei} we can
deduce the wildness of $\EO (g')$ from that of $\EO (g)$.
The two key ingredients of our work on blocks and quivers are
M\"{u}ller's theorem and skew group algebras. The former, which
was also fundamental to \cite{brogor1} and which is restated here
as Theorem \ref{muller}, implies that the blocks of $\EO (g)$ are
in bijection with the maximal ideals of $Z(\EO)$ lying over
$\mathfrak{m}_g$. The latter feature because the algebras $\EO
(g)$ are skew group algebras. This follows from the following, the
main result of Section 5, which is independent of the rest of the
paper and which may be of interest in other contexts.
\begin{itemize}
\item (Theorem \ref{skewgr}) Let $k$ be an algebraically closed
field and let $R$ be a finite dimensional $k$-algebra whose
irreducible modules are permuted simply transitively by a finite
abelian group $G$ of $k$-algebra automorphisms of $R$, with
$\mathrm{char} k$ coprime to $|G|$. Then $R$ is isomorphic to a
matrix algebra over a skew group algebra $S_1 *G$ with $S_1$
scalar local. The blocks and quiver of $R$ are determined by the
conjugation action of $G$ on the Jacobson radical of $S_1$.
\end{itemize}
\subsection{}
\label{I5} The paper also includes a series of results,
paralleling those in \ref{I3}, for the quantised Borel algebras
$\UM$ and $\UP$ at an $\ell$th root of unity $\epsilon$. Thus, see
Theorem \ref{tamewild} and Corollary \ref{utype} for their
representation type, and Theorems \ref{Upquiv} and \ref{ublocks}
for their blocks and quivers. The quantised analogues of the
enveloping algebras of the positive and negative Borel subalgebras
$\mathfrak{b}^+$ and $\mathfrak{b}^-$ of $\mathfrak{g} =
\mathrm{Lie}(G)$ are closely connected to $\EO$, thanks to the
dual of the multiplication map $m:B^+ \times B^- \longrightarrow
G$. Coupled with the Hopf self-duality of the quantised Borel
algebras, this yields an algebra embedding \cite[Section
6]{declyu1} of $\EO$ into $\UM \otimes \UP$. Corresponding to the
central subalgebra $\mathcal{O}[G]$ of $\EO$ are central
subalgebras $\mathcal{O}[B^+]$ and $\mathcal{O}[B^-]$ of $\UM$ and
$\UP$ respectively, and the above embedding specialises to the
finite dimensional factor algebras - see Proposition \ref{bigcell}
for details.
Despite the apparently simpler structure of the quantised
enveloping algebras of the Borels as compared with $\EO$, our
results for the former are in many cases weaker than for the
latter. There are at least two reasons for this: the coincidence
of Azumaya points with smooth points is not in general valid for
the centres of quantised enveloping algebras (the details are laid
out in Proposition \ref{azyouwere}); and  in general $Z(\UP)$ is
not known. (This situation has been rectified in \cite{gorDUR} where, in  particular, $Z(\UP)$ is described. Often, but not always, $Z(\UP) = \mathcal{O}[B^-]$ -
see Theorem \ref{ublocks}(ii).)
\subsection{}
\label{Iextra} The contents are  arranged as follows. In Section 2
notation is fixed and earlier work is recalled in the form most
useful for present purposes. We also provide a proof of the freeness of $\EO$ over its central subalgebra $\mathcal{O}[G]$. In Section 3 the centre of $\EO$ is
analysed with enough care to allow the proof of Theorem
\ref{azmor}, describing $\EO (g)$ when $g$ is fully Azumaya. In
Section 4 the representation type of $\EO (g)$ and of $\UP (b)$ is
determined. Section 5
contains the interlude on skew group algebras, culminating in
Theorem \ref{skewgr} and Proposition \ref{blquiv}. Sections 6 and
7 are concerned  with the block and quiver structure of
(respectively) $\UP (b)$ and $\EO (g)$. Both these sections
contain a number of examples.
The final three paragraphs of the Introduction suggest three
directions in which one might hope to extend the work described
here.
\subsection{}
\label{I6} The results of this paper show that the bundle of algebras $\{ \EO (g) : g \in G \}$ is
a partially ordered collection of successive degenerations,
progressing from the semisimple artinian algebras for $g \in
X_{w_0,w_0}$, the big cell, where
\[
\EO (g) \quad \cong \quad \left( \mathrm{Mat}_{\ell^N}(
\C)\right)^{\oplus(\ell^r)},
\]
towards the most degenerate algebras, for $g \in X_{e,e}$, where
the $\ell^r$ irreducible modules are one-dimensional and there is
only one block. This progressive degeneration is closely tied to
the Bruhat-Chevalley order on $W \times W$, (see Lemma
\ref{degen}). We exploit this perspective in analysing
representation type, for example, (as outlined in (\ref{I4})), but
it seems likely that more use can be made of similar arguments.
A similar philosophy applies to other classes of algebras whose
representation theory exhibits a geometric stratification, such as
the quantised enveloping algebras $\UE$ and the modular enveloping
algebras $\mathcal{U}(\g)$, for $\g$ semisimple; but the positive
evidence in these cases is more meagre than for the function
algebras, which - thanks to their structure as Galois coverings of
$\mathcal{O}[G]$ - provide a tractable testing ground for
techniques and conjectures to apply to the more difficult cases.
\subsection{}
\label{I7} A second  aspect where further work may prove fruitful
concerns the relations between the  structure of $\EO$ and
$\mathcal{O}_q [G]$, where $q$ is generic. The primitive ideals of
the generic algebras are also stratified by the orbits
$X_{w_1,w_2}$ of the double Bruhat cells, as shown in \cite{jos,
jos2}. Attempts have been made to determine the (second layer)
links between these primitive ideals, \cite{brogoo2, jos3}, but
the results remain incomplete. The work on blocks for $\EO$
presented here seems to support the belief that the links between
$\EO$-irreducibles are given by the ``lifts" modulo $\ell$ of
links
between primitive ideals in the corresponding stratum in
$\mathcal{O}_q [G]$. It may be that the results here can lead to
the formulation of the correct conjecture for the structure of the
link group in the generic case.
\subsection{}
\label{I8} Since the isomorphism type of $\EO (g)$ is determined
by the $T$-orbit $X_{w_1,w_2}$ of leaves to which $g$ belongs, it
is very natural to ask precisely what information regarding
$(w_1,w_2)$ suffices to determine the isomorphism type of $\EO
(g)$. It's already clear from Theorem \ref{simplefun}(b)(ii) that
$\ell(w_1),  \ell(w_2)$ and $s(w_2^{-1}w_1)$ are needed, but Corollary \ref{oblocks} indicates
that card$\mathfrak{S}(w_1,w_2)$ may be required also. And indeed
this is so - we show by example in \ref{oex}, with $G =
SL_{4}(\C )$, that card$\mathfrak{S}(w_1,w_2)$ isn't a function of
the other invariants listed above. Thus it remains an interesting
open problem to determine a ``minimal" set of isomorphism
invariants, in terms of Weyl group data, for the algebras $\EO
(g)$.
\section{Notations and Recollections}
\subsection{}
\label{not} Let $C=(a_{ij})$ be a Cartan matrix of finite type
having rank $r$ and let $(d_1, \ldots ,d_r)\in \mathbb{N}^r$ have
coprime entries such that $(d_ia_{ij})$ is symmetric. Let $\g$ be
the semisimple Lie algebra over $\mathbb{C}$ defined by $C$ and
let $\g = \mathfrak{n}^-\oplus \mathfrak{h}\oplus \mathfrak{n}^+$
be its triangular decomposition. Let $P$ and $Q$ be the weight and
root lattices of $\g$ and let $(\, , \,)$ be the associated
non-degenerate bilinear form. Let $\{ \alpha_1, \ldots
,\alpha_r\}$ be a set of simple roots determined by $C$ and let
$\{ \varpi_1,\ldots ,\varpi_r\}$ be the corresponding fundamental
weights of $P$. We have $(\varpi_i , \alpha_j) = \delta_{ij}d_i$.
Let $G$ be the simply-connected, semisimple algebraic group over
$\mathbb{C}$ associated with $C$. We have  Borel subgroups $B^+$
and $B^-$ of $G$ such that
$\text{Lie}(B^{\pm})=\mathfrak{n^{\pm}}\oplus \mathfrak{h}$. Let
$T=B^+\cap B^-$, a maximal torus of $G$. The Weyl group of $G$
(with respect to $T$) is $N_G(T)/T$. This can be identified with
the Weyl group associated with $C$. The Weyl group acts on both
$P$ and $Q$ and the form $(\, , \,)$ is $W$-invariant.
There is a stratification of $G$:
\[
G= \coprod_{w_1,w_2\in W} X_{w_1,w_2}
\]
where $X_{w_1,w_2}=B^+w_1B^+\cap B^-w_2B^-$. This restricts to a
stratification of $B^-$:
\[
B^- = \coprod_{w\in W} X_{w, e}.
\]
Any element $w\in W$ can be written as a product of simple
reflections or as a product of (arbitrary) reflections. We let
$\ell (w)$ (respectively $s(w)$) equal the minimal length of an
expression for $w$ as a product of simple (respectively arbitrary)
reflections. The longest word with respect to $\ell$ will be
denoted $w_0$; recall that $\ell(w_0) = N$, where $N =
\mathrm{dim}_{\C}(\mathfrak{n}^+)$, the number of positive roots.
The function $s:W\longrightarrow \mathbb{N}$ is called the rank
function. It coincides with the codimension of
$\mathbb{Q}\otimes_{\mathbb{Z}}P^w$ in
$\mathbb{Q}\otimes_{\mathbb{Z}}P$, where we write $P^w$ to denote
the elements of $P$ fixed by $w$.
Let $h$ be the Coxeter number of $W$. This equals the order in $W$
of the product of the simple reflections.
Throughout this paper $\ep \in \mathbb{C}$ will be a primitive
$\ell$th root of unity for some natural number $\ell >1$. Let
$\theta = \sum a_i\alpha_i$ be the highest root of $\g$.  We will
always require that $\ell $ is good, that is $\ell$ is odd and
prime to the integers $a_i$ and $d_i$ for $1\leq i\leq r$.
\subsection{}
\label{free}
Let $\UP$ (respectively $\UM$) be the non-negative (respectively
non-positive) subalgebra of the (simply-connected)
quantised enveloping algebra at a root of unity, $\ep$, associated
to $C$, as defined in \cite{conkacpro59}. Let $\EO$ be the quantised function
algebra at a root of unity, $\ep$, associated to $C$, as defined
in \cite{declyu1}.
The ring of regular functions on $B^-$, $\mathcal{O}[B^-]$, is a
central sub-Hopf algebra of $\UP$. Moreover, $\UP$ is free as a module
over $\mathcal{O}[B^-]$ of
rank $\ell^{\text{dim} B^-}$, \cite[Corollary 3.3(b)]{deck}.
Similarly, the ring of regular
functions on $G$, $\mathcal{O}[G]$, is a central sub-Hopf algebra
of $\EO$, \cite[Theorem 6.4]{declyu1}. We'll use $\theta$
to denote the embedding of $\mathcal{O}[G]$ into $Z(\EO)$, and for
the embedding of $\mathcal{O}[B^-]$ into $Z(\UP)$. In both cases
we'll denote the induced map $\theta^*$ on maximal spectra by
$\pi$.
\begin{prop}
As an $\mathcal{O}[G]$-module, $\EO$ is free of rank $\ell^{\dim G}$.
\end{prop}
\begin{proof}
Thanks to \cite[Theorem 7.2]{declyu1} $\EO$ is a projective
$\mathcal{O}[G]$-module of rank $\ell^{\dim G}$. By \cite{marlin}
the Grothendieck group of projective modules over $\mathcal{O}[G]$
is trivial, in other words
\[
K_0(\mathcal{O}[G]) \cong \mathbb{Z}.
\]
In particular, if $P$ is a projective $\mathcal{O}[G]$-module
whose rank is greater than the Krull dimension of $G$, then $P$ is
necessarily free, \cite[Theorem 11.3.7]{McC-Rob}. Since $\ell >1$
we have $\text{Kdim}\mathcal{O}[G] = \dim G < \ell^{dim G} =
\text{rank}\EO$, so the proposition follows.
\end{proof}
\subsection{}
\label{simplefun}
Let $b\in B^-$ and $g\in G$ and let $\mathfrak{m}_b\triangleleft
\mathcal{O}[B^-]$ and $\mathfrak{m}_g\triangleleft \mathcal{O}[G]$ be
the maximal ideals associated to these points. We define
\[
\UP (b) \equiv \frac{\UP}{\mathfrak{m}_b \UP}, \qquad \EO(g) \equiv
\frac{\EO}{\mathfrak{m}_g \EO}.
\]
In view of \ref{free}, these algebras have
$\mathbb{C}$-dimension $\ell^{\text{dim}B^-}=\ell^{N+r}$ and
$\ell^{\text{dim}G}=\ell^{2N+r}$ respectively.
\begin{thm}\cite[Theorem 4.4]{conkacpro59},\cite[Section
9]{declyu1},\cite[Theorem 4.4 and Proposition 4.10]{decpro49}
Let $b, b'\in X_{w,e}$ and $g,g'\in X_{w_1,w_2}$ for some
$w,w_1,w_2\in W$.\\
(a)(i) There is an algebra isomorphism $\UP (b)\cong \UP (b')$. \\
(ii) There are precisely $\ell^{r-s(w)}$ simple $\UP (b)$-modules and
each simple module has dimension $\ell^{\frac{1}{2}(\ell(w)+s(w))}$.\\
(b)(i) There is an algebra isomorphism $\EO (g)\cong \EO (g')$. \\
(ii) There are precisely $\ell^{r-s(w_2^{-1}w_1)}$ simple $\EO
(g)$-modules and each simple module has dimension $\ell^{\frac{1}{2}(\ell(w_1) +
\ell(w_2) +s(w_2^{-1}w_1))}$.
\end{thm}
\subsection{}
\label{centre}
We recall the description of the centre $Z(\EO)$ of $\EO$ given in
\cite{enr} and
\cite[Appendix]{declyu1}. Let $U_q$ be the quantised enveloping algebra associated with
Cartan matrix $C$, defined over $\mathbb{C}(q)$, with $q$ an
indeterminate. For $1\leq i\leq r$ let $L(\varpi_i)$ be the simple
$U_q$-module of type 1 with highest weight $\varpi_i$. Let
$v_{\varpi_i}$ (respectively $f_{-w_0\varpi_i}$) denote the highest
weight vector of $L(\varpi_i)$ (respectively $L(\varpi_i)^*$) and let
$v_{-\varpi_i}$ (respectively $f_{w_0\varpi_i}$) denote the lowest
weight vector of $L(-w_0\varpi_i)$ (respectively
$L(-w_0\varpi_i)^*$). These are well-defined up to scalar
multiplication. We define the (quantum) matrix coefficients
\[
\bi = c^{\varpi_i}_{f_{-w_0\varpi_i},v_{\varpi_i}}, \qquad
\ci=c^{-w_0\varpi_i}_{f_{w_0\varpi_i},v_{-\varpi_i}}.
\]
These elements can (and will) be considered as elements of $\EO$ after specialisation of an appropriate integral form. Let $Z_q$ be the subalgebra of $\EO$ generated by the elements $b_i^kc_i^{\ell - k}$ for $1\leq i \leq r$ and $0\leq k \leq \ell$.
\begin{thm}[Enriquez]
The centre $Z(\EO )$ of $\EO$ is isomorphic to $\mathcal{O}[G]\otimes_{\mathcal{O}[G]\cap Z_q}Z_q$.
\end{thm}
\begin{rems}
(1) As a $\C$-algebra $\mathcal{O}[G]\cap Z_q$ is generated by $\{
b_{\varpi_i}^{\ell}, c_{\varpi_i}^{\ell} : 1\leq i \leq r\}$, and
is isomorphic to $\C [X_1,\ldots X_r, Y_1,\ldots ,Y_r ]$.
\newline
\indent (2) We can identify $Z_q$ with $\C[\alpha_i(k) :1\leq
i\leq r,0\leq k\leq \ell]/I$ where $I$ is generated by the
elements, for $1\leq i\leq r$,
\begin{equation*}
\begin{cases}
\alpha_i(k)\alpha_i(k')-\alpha_i(0)\alpha_i(k+k') &\text{if }
k+k'\leq \ell, \\
\alpha_i(k)\alpha_i(k')-\alpha_i(\ell)\alpha_i(k+k'-\ell)
&\text{if } k+k'>\ell,
\end{cases}
\end{equation*}
and the identification maps $\alpha_{i}(k)$ to
$b_{i}^{k}c_{i}^{\ell - k}$, so that $\mathcal{O}[G]\cap
Z_q\hookrightarrow Z_q$ corresponds to $X_i \longmapsto
\alpha_i(0)$ and $Y_i\longmapsto \alpha_i(\ell )$.
\newline
\indent (3) $Z_q\cong \bigotimes_{i=1}^r Z_q^i$, where $Z_q^i$ is
the algebra generated by $\{\alpha_i(k):0\leq k\leq \ell\}.$
\newline
\indent (4) As an $\mathcal{O}[G]\cap Z_q$-module, $Z_q$ is
finitely generated and free: indeed each $Z_q^i$ is free over
$\mathcal{O}[G]\cap Z_q^i$ with basis $\{1\} \cup \{\alpha_i(k) :
1\leq k\leq \ell -1 \}$.
\end{rems}
\subsection{}
\label{maxfactor}
Given any simple $\EO$-module $V$ the centre of $\EO$ acts by scalar
multiplication thanks to Schur's lemma. Thus there is an algebra map
\[
\zeta_V :Z(\EO) \longrightarrow \mathbb{C}
\]
which we call the central character of $V$. We define the
following subset of $\text{Maxspec}(Z(\EO))$:
\[
\mathcal{A}_{\EO} \equiv \{ \text{ker}(\zeta_V) : V \textit{ is a
simple $\EO$-module of maximal dimension} \}.
\]
This set is non-empty and open in $\text{Maxspec}(Z(\EO))$ and is
called the \textit{Azumaya locus of} $\EO$. Thanks to Theorem
\ref{simplefun} we have
\begin{eqnarray*}
\label{azum}
 \mathcal{A}_{\EO} = \{ \text{ker}(\zeta_V): &&V
\textit{ is a simple $\EO (g)$-module for $g\in X_{w_1,w_2}$} \\
&& \textit{with $\ell(w_1)+\ell(w_2)+s(w_2^{-1}w_1) = 2N$}\}.
\end{eqnarray*}
Following \cite[Section 2.5]{brogor1} we also define the
\textit{fully Azumaya locus} $\mathcal{F}_{\EO}$ of
$\mathcal{O}[G]$ with respect to $\EO$ to consist of those maximal
ideals $\mathfrak{m}_g$ of $\mathcal{O}[G]$ (or equivalently those
elements $g$ of $G$) such that every irreducible $\EO (g)$-module
has dimension $\ell^N$. In the notation of \ref{simplefun},
\[
\pi(\mathcal{A}_{\EO}) = \mathcal{F}_{\EO} = \{ g \in G : g \in
X_{w_1,w_2}, \textit{ with } \ell(w_1) + \ell(w_2) + s(w_2^{-1} w_1)
= 2N \},
\]
and similarly for $\UP$.
\begin{thm}\cite[Corollary 2.7]{brogor1}
Let $g\in\mathcal{F}_{\EO}$. Then
\[
\EO (g) \cong \mathrm{Mat}_{\ell^N}\left(
\frac{Z(\EO)}{\mathfrak{m}_gZ(\EO)}\right).
\]
\end{thm}
\subsection{}
\label{az}
There is an alternative description of $\mathcal{A}_{\EO}$.
\begin{thm}\cite[Theorem C]{brogoo}
The Azumaya locus $\mathcal{A}_{\EO}$ of $\EO$ coincides with the
smooth locus of $\mathrm{Maxspec}(Z(\EO))$.
\end{thm}
In Section 3 we will use Theorems \ref{maxfactor} and \ref{az} to
describe the algebras $\EO (g)$ for all $g$ in
$\mathcal{F}_{\EO}$.
\subsection{}
\label{azyouwere}
 Theorem \ref{az} is not in general true for
$\UP$.
\begin{prop}
The Azumaya locus
$\mathcal{A}_{\UP}$ of $\UP$ coincides with smooth locus of
$\mathrm{Maxspec}(Z(\UP))$ if and only if $C$ is a Cartan matrix of
type $A_{2n}$.
\end{prop}
\begin{proof}
Let $Z = Z(\UP)$. In these circumstances \cite[Theorem 3.8]{brogoo} states that the Azumaya locus $\mathcal{A}_{\UP}$ of $\UP$
coincides with the smooth locus of $\mathrm{Maxspec}(Z)$ if $\UP$ is Azumaya in codimension one, that is the set of points of
$\mathrm{Maxspec}(Z)$ which are not annihilators of simple
$\UP$-modules of maximal dimension has codimension at least two, see
\cite[Corollary 1.8]{brogoo}.
Since $Z$ is the centre of a maximal order it is integrally
closed, \cite[Theorem 5.1.10(b)]{McC-Rob}. In particular
$\mathrm{Maxspec}(Z)$ is smooth in codimension one. This means
that the converse of \cite[Theorem 3.8]{brogoo} is also true:
if $\mathcal{A}_{\UP}$ coincides with the smooth locus of
$\mathrm{Maxspec}(Z)$ then $\UP$ is necessarily Azumaya in
codimension one.
The map induced by inclusion
\[
\pi : \mathrm{Maxspec}(Z) \longrightarrow \mathrm{Maxspec}(Z_0) = B^-,
\]
is surjective with finite fibres. The simple $\UP$-modules lying
over $b\in X_{w,e}$ all have dimension $\ell^{\frac{1}{2}(\ell(w)
+ s(w))}$ and the maximal dimension of a simple $\UP$-module is
$\ell^{\frac{1}{2}(N + s(w_0))}$. By \cite[Theorem 1.1]{FZ} the
variety $X_{w,e}$ has dimension $\ell(w) + r$, so
$\mathrm{Maxspec}(Z)$ is stratified by pieces $\pi^{-1}(X_{w,e})$
of dimension $\ell(w) +r$ over which the representation theory is
constant. Hence $\UP$ is Azumaya in codimension one if and only if
$\ell(w) + s(w) = N + s(w_0)$ for all $w\in W$ such that $N -
\ell(w) \leq 1$. In other words we need only check that $s(w_0s_i)
= s(w_0) + 1$ for all $1\leq i\leq r$. By \cite[Lemma 7.6]{gor4}, however,
this is equivalent to the condition $-w_0 (\alpha_i) \neq
\alpha_i$ for all $1\leq i\leq r$, in other words the involution
$-w_0$ has no fixed points. This happens if and only if $C$ is a
Cartan matrix of type $A_{2n}$.
\end{proof}
\subsection{}
\label{muller}
The following result was proved in a ring-theoretic setting by
M\"{u}ller, \cite{mul}. A discussion of the form given here can be
found in \cite[Paragraph 2.10]{brogor1}.
\begin{thm}
Let $b\in B^-$ and $g\in G$. The blocks of $\UP (b)$ are in natural
correspondence with the maximal ideals of $Z(\UP)$ lying over
$\mathfrak{m}_b$. Similarly, the blocks of $\EO (g)$ are in
correspondence with the maximal ideals of $Z(\EO)$ lying over
$\mathfrak{m}_g$.
\end{thm}
\subsection{}
\label{tamewild}
We recall the notion of representation type of a finite dimensional
algebra $T$:
\\
\indent (i) $T$ has finite representation type if there are a
finite number of mutually non-isomorphic indecomposable
$T$-modules;\\ \indent (ii) $T$ has tame representation type if
$T$ does not have finite representation type and if, for each
dimension $d> 0$, there is a finite collection of
$T-\mathbb{C}[x]$-bimodules $M_i$ which are free as right
$\mathbb{C}[x]$-modules such that every indecomposable $T$-module
of dimension $d$ is isomorphic to $M_i\otimes_{\mathbb{C}[x]} N$
for some $i$ and some simple $\mathbb{C}[x]$-module $N$;
\\
\indent (iii) $T$ has wild representation type if there is a
finitely generated $T-\mathbb{C}<x,y>$-bimodule $M$ which is free
as a right $\mathbb{C}<x,y>$-module such that the functor
$F(N)=M\otimes_{\mathbb{C}<x,y>}N$ from the category of finite
dimensional $\mathbb{C}<x,y>$-modules to the category of finite
dimensional $T$-modules preserves indecomposability and
isomorphism classes.
By \cite{dro} $T$ falls into precisely one of the above classes: we will say
$T$ is finite, tame or wild as appropriate.
\begin{thm}\cite[Theorem 7.1]{gor4} In addition to the usual hypotheses on $\ell$, assume that $\ell >
h$. Let $b\in X_{w,e}$ and $g\in X_{w_1,w_2}$ for some $w,w_1,w_2
\in W$.\\ (a)(i) If $\ell(w) > N-2$ then $\UP (b)$ has finite
representation type.\\ (ii) If $\ell (w) < N - 2$ then $\UP (b)$
has wild representation type.\\ (b)(i) If
$\ell(w_1)+\ell(w_2)>2N-2$ then $\EO (g)$ has finite
representation type. \\ (ii) If $\ell (w_1)+\ell(w_2) <2N-2$ then
$\EO (g)$ has wild representation type.
\end{thm}
In Theorem \ref{reptypethm2} and Corollary \ref{utype} we'll
complete the determination of the representation type of $\EO (g)$
and $\UP (b)$ and remove the restriction $\ell > h$.
\subsection{}
\label{description}
We have a description of the algebras occurring in Theorem
\ref{tamewild}(a)(i).
\begin{prop}\cite[Theorem 7.7]{gor4}
Assume $\ell >h$. Let $b\in X_{w,e}$.\\
(a) Assume $\ell(w)=N$. Then $w=w_0$ and there is an algebra isomorphism
\[
\UP (b) \cong \bigoplus_{j=1}^{\ell^{r-s(w_0)}}
\mathrm{Mat}_{\ell^{N+s(w_0)}}\left( \mathbb{C}\right) .
\]
(b) Assume $\ell(w) = N-1$. Then $w=w_0s_i$ for some $i$. There are
two cases:\\
(i) $w_0(\alpha_i) = -\alpha_i$ - there is an algebra isomorphism
\[
\UP (b) \cong \bigoplus_{j=1}^{\ell^{r-s(w_0)}}
\mathrm{Mat}_{\ell^{\frac{1}{2}(N+s(w_0)-2)}}\left(\BARUP
(\mathfrak{sl}_2)\right);
\]
(ii) $w_0(\alpha_i) \neq -\alpha_i $ - there is an algebra
isomorphism
\[
\UP (b) \cong
\bigoplus_{j=1}^{\ell^{r-s(w_0)-1}}\mathrm{Mat}_{\ell^{\frac{1}{2}(N+s(w_0))}}\left(\frac{\mathbb{C}[X]}{(X^{\ell})}\right).
\]
\end{prop}
We will see in Lemma \ref{lessres} that the restriction $\ell>h$ can be removed
from this proposition. There is also a corresponding description for the algebras in Theorem
\ref{tamewild}(b)(i) given in \cite[Theorem 7.4]{gor4}. One can recover this from Theorem \ref{azmor}.
\subsection{}
\label{appnot}
Let us recall the definition of the (right) winding automorphisms for $\UP (b)$
(respectively $\EO (g)$). Given a character of $\UP$ factoring through
$\BARUP$,
\[
\chi : \UP \longrightarrow \BARUP \longrightarrow \mathbb{C},
\]
we define an automorphism $\tau_{\chi}$ of $\UP$ by
\[
\tau_{\chi}(x) = \sum_{(x)} x_{(1)}\chi( x_{(2)}).
\]
Since $Z_0 \subseteq \UP$ is a sub-Hopf algebra on which $\chi$
agrees with the augmentation $\ep$, we see that $\tau_{\chi}$ acts
as the identity on $Z_0$. Therefore $\tau_{\chi}$ induces an
automorphism on $\UP (b)$ for any $b \in B^-$. It is
straightforward to check that for any $\UP (b)$-module $M$ the
twisted module $^{\tau_{\chi}}M$ is isomorphic to $M \otimes
\C_{\chi}$. One argues similarly for $\EO$.
The characters of $\BARUP$ are parametrised by $Q_{\ell} =
Q/\ell Q$. Namely, for any element $\mu\in Q$ we have the one dimensional
representation given by
\[
E_i . 1 = 0 \quad \textrm{and} \quad K_{\lambda} . 1
=\ep^{(\lambda , \mu)}.
\]
These representations are different for different choices of coset
representative of $\ell Q$ in $Q$, and every irreducible
$\BARUP$-module arises in this way by Theorem \ref{simplefun}(a)(ii) applied in the case
$w=e$. More generally,
the following theorem is proved in \cite[Theorems 4.5 and 4.10]{decpro49}.
\begin{thm}
(i) Let $b\in X_{w,e}$. If $\mu \in Q$ and $S$ is a simple $\UP (b)$-module, then
$^{\tau_{\mu}}S \cong S$ if and only if $(\lambda , \mu )\in
\ell\mathbb{Z}$ for all $\lambda \in P^w$.
\\
\indent (ii) Let $g\in X_{w_1,w_2}$. If $\mu \in Q$ and $S$ is a simple $\EO (g)$-module, then $^{\tau_{\mu}}S
\cong S$ if and only if $(w_1(\lambda), \mu)\in \ell\mathbb{Z}$
for all $\lambda \in P^{w_2^{-1}w_1}$.
\end{thm}
\subsection{}
\label{ellform}
We show
that a subgroup of $Q_{\ell}$ acts simply transitively on the simple
$\UP (b)$-modules (respectively simple $\EO(g)$-modules).
\begin{lem}
(i) The
nondegenerate form
\begin{eqnarray}
P \times Q \longrightarrow \mathbb{C}(q) : (\alpha,
\beta)\longmapsto q^{(\alpha,\beta)}
\end{eqnarray}
induces a nondegenerate form
\begin{eqnarray}
P_{\ell} \times Q_{\ell} \longrightarrow \mathbb{C} : (\alpha,
\beta)\longmapsto \epsilon^{(\alpha,\beta)}.
\end{eqnarray}
\indent (ii) There is an elementary abelian $\ell$-subgroup of $Q_{\ell}$ which acts simply transitively on the
simple $\UP (b)$-modules. If $\ell$ is prime to the order of $w$ the
subgroup can be chosen to be $Q^w/\ell Q^w$. \\
\indent (iii) There is an elementary abelian  $\ell$-subgroup of $Q_{\ell}$
which acts simply transitively on the simple $\EO (g)$-modules. If $\ell$ is
prime to the order of $w_2^{-1}w_1$ then this subgroup can be chosen
to be $Q^{w_2w_1^{-1}}/\ell Q^{w_2w_1^{-1}}$.\end{lem}
\begin{proof}
(i) We overline to indicate images modulo $\ell$ in
a $\mathbb{Z}$-module. Now $P_{\ell} =
\sum_{j}\overline{\mathbb{Z}\varpi_{i}}$ and $Q_\ell
= \sum_{j}\overline{\mathbb{Z}\alpha_{i}}$, with
$(\overline{\alpha_i},\overline{\varpi_j}) = \delta_{ij}d_i$,
 a unit in $\overline{\mathbb{Z}}$ when $i = j$. Thus $P_{\ell} = (Q_{\ell})^*$ and $Q_{\ell} = (P_{\ell})^*$, as
 claimed.\\
\indent (ii) That $Q_{\ell}$ acts transitively on the simple $\UP
(b)$-modules follows from \cite[Theorem 4.5]{decpro49}. Let $b \in X_{w,e}$, and write $\mathbb{Z}'$ for $\mathbb{Z}[d_1^{-1}, \ldots , d_r^{-1}]$. We write $M' =
M\otimes_{\mathbb{Z}}\mathbb{Z}'$ for a $\mathbb{Z}$-module $M$. Since $(\alpha_i,\varpi_{j}) = d_{i}\delta_{ij}$, there is a perfect
pairing
\[
(\, , \,) : P'\times Q' \longrightarrow \mathbb{Z}'.
\]
Suppose first that $P'=P_1'\oplus P_2'$. Since $Q' = \textrm{Hom}_{\mathbb{Z}'}(P', \mathbb{Z}')$ via the perfect pairing, $Q' = P_{1}'{}^{\perp} \oplus
P_{2}'{}^{\perp}$.\\
\indent Let $P'{}^w \subseteq P'$ and note that $P'/P'{}^w$ is
torsion-free since
$n\lambda \in P'{}^w$ implies that $\lambda \in P'{}^w$. Hence $P' = P'{}^w \oplus P_2'$ for some
$P_2'$. Thus $Q' = Q_1 \oplus Q_2$, where $Q_1 = (P'{}^{w}){}^{\perp}$
and $Q_2 = (P_2'){}^{\perp}$. Let $\mu \in Q$ be such that $(P^w, \mu)
\subseteq \ell\mathbb{Z}$.  Since the integers $d_i$ are prime to
$\ell$ (thanks to our continuing hypothesis on $\ell$ given in
Paragraph \ref{not})
this is equivalent to $(P'{}^w , \mu) \subseteq \ell\mathbb{Z}'$. Write $\mu = \mu_1 + \mu_2$ with $\mu_i \in Q_i$. Then for $\lambda \in P'{}^w, (\lambda, \mu_2 ) = (\lambda , \mu) \in \ell
\mathbb{Z}'$, so that $(P', \mu_2) \subseteq
\ell\mathbb{Z}'$. Since the pairing is perfect, this forces  $\mu_2 \in \ell Q_2$.
Hence $Q_2/\ell Q_2$ operates simply transitively on the
irreducible $\UP (b)$-modules by Lemma \ref{appnot}(i).
Now suppose that the order of $w$ is prime to $\ell$. Let
$\mathbb{Z}'' = \mathbb{Z}'[\text{ord}(w)^{-1}]$. Working over $\mathbb{Z}''$, we can argue as
 above to find a
decomposition $P'' = P''{}^w \oplus P_2''$ with $P_2''$
$<w>$-invariant. We claim
that $Q_2'' \subseteq Q''{}^w$. This follows from the observation
that
$(\mu - w\mu , P''{}^w) = 0 = (\mu - w\mu , P_2'')$, the first equality
by $<w>$-invariance of the elements of $P''{}^w$, the second by
orthogonality and the $<w>$-invariance of $P_2''$. Thus we have a
factorisation
\[
Q_2 /\ell Q_2 \twoheadrightarrow Q_2/Q_2\cap \ell Q^w \hookrightarrow
Q^w /\ell Q^w.
\]
Since the right hand side and the left hand side both have $\ell
^{r-s(w)}$ elements this completes the proof of (ii).\\ \indent
(iii) is proved entirely similarly.
\end{proof}
\subsection{}
Recall that if $S$ is a finite dimensional algebra on which a
group $G$ acts by algebra automorphisms then we can form the skew
group algebra of $S$ by $G$, written $S\ast G$. As a left
$S$-module this is free with basis $g\in G$. Multiplication is
given by extension of the formula $s^{g} g = g s$, for $s \in S$
and $g\in G$ and where $s^{g}$ denotes the action of $g$ on $s$.
\subsection{}
\label{bigcell} Let $m:B^+\times B^- \longrightarrow G$ be the
multiplication map. Then $m$ is a principal $T$-bundle onto the
open, dense subset $B^+B^-$ of $G$. Let $b_1\in X_{w_1,e}$ and
$b_2\in X_{e,w_2}$ be unipotent, (the double Bruhat cells are
$T$-invariant, so that we can find such representatives), and let
$g=m(b_2,b_1)\in X_{w_1,w_2}$.
\begin{prop}
\cite[Section 2.9]{gor4} Let $w_1, w_2\in W$, $g\in G$, $b_1\in B^-$
and $b_2 \in B^+$ be as above. Then there is an algebra
isomorphism
\[
\EO (g)\ast \mathbb{Z}_{\ell}^r \longrightarrow \UM (b_2)\otimes \UP
(b_1).
\]
\end{prop}
\section{The Azumaya locus of $\EO$}
\subsection{}
Recall $\mathcal{O}[G]$ and $Z_q$, the central subalgebras of $\EO$
introduced in \ref{simplefun} and \ref{centre} respectively. We will denote $\text{Spec}Z_q$ by $X$ and $\text{Spec}\mathcal{O}[G]\cap Z_q$ by $U$. By Remark \ref{centre}(1) $U$ is isomorphic to $\mathbb{A}^{2r}$.
\begin{lem}
\label{sing1}Let $\pi :X\longrightarrow U$ be the (surjective) morphism induced by the inclusion $\mathcal{O}[G]\cap Z_q\longrightarrow Z_q$. Write the points of $U$ as $2r$-tuples $(b_1,\ldots ,b_r,c_1,\ldots c_r)$ under the identification in Remark \ref{centre}(1). If $p\in U$ is such that $b_i
=0=c_i$ for some $i$ then every point of $\pi ^{-1}(p)$ is singular in $X$.
\end{lem}
\begin{proof}
It is clear that $\pi$ factorises as $\pi_1\times \ldots \times \pi_r$ where $\pi_i: \text{Spec}Z_q^i \longrightarrow \text{Spec}\mathcal{O}[G]\cap Z_q^i$.
 It is therefore enough to prove this for the case $r=1$.
By Remark \ref{centre}(2) we can consider $X$ as an affine variety embedded in
 $\mathbb{A}^{\ell +1}$ (with co-ordinate functions $\alpha_1(k)$). Under this identification $\pi$ takes $(a_0,\ldots ,a_{\ell})$ to $(a_0,a_{\ell})$.
Note that $\pi^{-1}((0,0))= (0,\ldots ,0)$. Indeed the equations
\[
\alpha_i(k)\alpha_i(k')=
\begin{cases}
\alpha_i(0)\alpha_i(k+k') &\text {if }k+k'\leq \ell, \\
\alpha_i(\ell )\alpha_i(k+k'-\ell ) &\text {if} k+k'> \ell,
\end{cases}
\]
show that $a_k^2=0$ as required.
Define $f_{k,k'}$ as follows:
\[
f_{k,k'} =
\begin{cases}
\alpha_i(k)\alpha_i(k') - \alpha_i(0)\alpha_i(k+k') &\text {if }k+k'\leq \ell, \\
\alpha_i(k)\alpha_i(k') - \alpha_i(\ell )\alpha_i(k+k'-\ell ) &\text {if } k+k'> \ell.
\end{cases}
\]
Recall that for $p= (a_0,\ldots ,a_{\ell})\in X$ we define
\[
f_{k,k',p}^{(1)} = \sum_{j=0}^{\ell} \frac{\partial f_{k,k'}}{\partial \alpha_1(j)}(p)(\alpha_1(j)-a_j).
\]
Then, by definition, the tangent space of $X$ at $p$ is
\[
T_pX = \cap_{1\leq k,k'\leq \ell -1}\{x \in \mathbb{A}^{\ell + 1}
: f_{k,k',p}^{(1)} = 0) .
\]
Since $f_{k,k'}$ is homogeneous of degree two it follows that for all $k,k'$
\[
f_{k,k',(0,0)}^{(1)}\equiv 0,
\]
which implies that $\dim (T_0X)=\ell +1$. Now $Z_q$ is finite over
$\mathcal{O}[G]\cap Z_q$, so that $X$ has dimension 2. Therefore
$(0,\ldots ,0)=\pi^{-1}((0,0))$ is a singular point.
\end{proof}
\subsection{}
The following proposition allows us to ignore some unfavourable points
of $X$.
\begin{prop}
\label{sing2}
Let $g\in G$ be such that $b_{\varpi_i}^{\ell} (g)=0 = c_{\varpi_i}^{\ell}(g)$ for some $i$, $1\leq i\leq r$. Then $\EO (g)$ is not Azumaya.
\end{prop}
\begin{proof}
Using Theorem \ref{az} it is enough to show that any maximal ideal of
$Z$ lying over $\mathfrak{m}_g$ is singular. Denote
$\text{Spec}(Z(\EO))$ by $Y$, the fibre product $X\times_U G$.
\newline
\noindent
\textbf{Claim. }
Let $\tilde{\pi} : Y\longrightarrow X$ be the projection map. If $x\in X$ is singular then any point of $\tilde{\pi} ^{-1}(x)$ is singular in $Y$.
\newline
\noindent
\textit{Proof of claim. }
For ease of notation let $R=Z_q$, $S=Z_q\cap \mathcal{O}[G]$ and $T=\mathcal{O}[G]$. Thus we are considering $R\otimes_S T$ where:
\newline
(i) as an $S$-module $R$ is finitely generated and free;
\newline
(ii) the algebra $S$ is smooth;
\newline
(iii) all algebras are affine domains.
Let $\mathfrak{m}_R \triangleleft R$ be the maximal ideal corresponding to $x\in X$, and suppose $M$ lies over $\mathfrak{m}_R$. Define $\mathfrak{m}_S=M\cap S$ and $\mathfrak{m}_T=M\cap T$, maximal ideals of $S$ and $T$ respectively.
We first show that
\begin{equation}
\label{eq:localise}
(R\otimes_S T)_M = R_{\mathfrak{m}_R}\otimes _{S_{\mathfrak{m}_S}}T_{\mathfrak{m}_T}.
\end{equation}
Note that $M = \mathfrak{m}_R\otimes_S T + R\otimes_S \mathfrak{m}_T$. This means in particular that if $y\in R \setminus\mathfrak{m}_R$ and $z\in T\setminus\mathfrak{m}_T$ then $y\otimes z\in R\otimes_S T\setminus M$. So we have an embedding
\begin{equation}
\label{eq:embedloc}
A:=R_{\mathfrak{m}_R}\otimes _{S_{\mathfrak{m}_S}}T_{\mathfrak{m}_T} \longrightarrow  (R\otimes_S T)_M.
\end{equation}
Let $x\in R\otimes_S T\setminus M$. If $A$ were local then $x\in
A\setminus MA$ would be an invertible element, since $MA$ is
maximal. Therefore the map (\ref{eq:embedloc}) would be an
isomorphism, proving the claim. Thus it is enough to show that $A$
is a local ring. Observe that by (i) the algebra
$R/\mathfrak{m}_SR$ is finite dimensional. Therefore localising
this at $\mathfrak{m}_R$ yields another finite dimensional algebra
$R_{\mathfrak{m}_R}/\mathfrak{m}_SR_{\mathfrak{m}_R}$. Nakayama's
lemma implies that $R_{\mathfrak{m}_R}$ is finite over
$S_{\mathfrak {m}_S}$.
Let $\overline{\mathfrak{m}}=\mathfrak{m}_TA$. Then
\[
A/\overline{\mathfrak{m}} \cong R_{\mathfrak{m}_R}/\mathfrak{m}_SR_{\mathfrak{m}_R}
\]
is finite dimensional, so Nakayama's lemma also implies that $A$ is finite over $T_{\mathfrak{m}_T}$. Therefore $\overline{\mathfrak{m}}$ is contained in the Jacobson radical of $A$. However, since $A/\overline{\mathfrak{m}}$ is local it follows that $A/Jac(A)$ is local, as required.
So we have proved (\ref{eq:localise}).
To complete the claim we must show that $A=R_{\mathfrak{m}_R}\otimes _{S_{\mathfrak{m}_S}}T_{\mathfrak{m}_T}$ has infinite global dimension. By hypothesis $R_{\mathfrak{m}_R}$ has.
We have a change of rings spectral sequence
\begin{equation}
\label{eq:corspecseq}
\x{A}{p}{k_A}{\x{R_{\mathfrak{m}_R}}{q}{A}{M}} \Longrightarrow \x{R_{\mathfrak{m}_R}}{p+q}{k_{\mathfrak{m}_R}}{M},
\end{equation}
for any $R_{\mathfrak{m}_R}$-module $M$.
By Frobenius reciprocity we have
\[
\x{R_{\mathfrak{m}_R}}{q}{A}{M}= \x{R_{\mathfrak{m}_R}}{q}{R_{\mathfrak{m}_R}\otimes _{S_{\mathfrak{m}_S}}T_{\mathfrak{m}_T}}{M} \cong \x{S_{\mathfrak{m}_S}}{q}{T_{\mathfrak{m}_T}}{M}.
\]
As $S_{\mathfrak{m}_S}$ is smooth there exists a natural number $Q$ such that
\[
\x{S_{\mathfrak{m}_S}}{q}{T_{\mathfrak{m}_T}}{M}= 0,
\]
for all $q > Q$ and all $S_{\mathfrak{m}_S}$-modules $M$.
If $A$ were smooth there would be a natural number $P$ such that
\[
\x{A}{p}{k_A}{M} = 0,
\]
for all $p >P$ and all $A$-modules $M$. Then by (\ref{eq:corspecseq}) we have
\[
\x{R_{\mathfrak{m}_R}}{n}{k_{\mathfrak{m}_R}}{M} =0,
\]
for all $n > P+Q$ and $R_{\mathfrak{m}_R}$-modules $M$, contradicting the singularity of $R_{\mathfrak{m}_R}$. Thus $A$ must be singular, proving the claim.
The proposition now follows from Lemma \ref{sing1} combined with the claim.
\end{proof}
\subsection{}
\label{allaz}
Now we can describe the algebras lying over the fully Azumaya
locus - that is, we describe $\EO (g)$ for $g \in
\mathcal{F}_{\EO}$, in the notation of \ref{az}.
\begin{thm}
\label{azmor} Suppose $g\in X_{w_1,w_2}\cap \mathcal{F}_{\EO}$; that is, $\ell(w_1) + \ell(w_2) + s(w_{2}^{-1}w_1) = 2N$.
Let $s=s(w_2^{-1}w_1)$. Then there is an algebra isomorphism
\[
\EO (g) \cong \bigoplus_{1}^{\ell ^{r-s}} \mathrm{Mat}_{\ell ^N}
\left( \frac{k[X_1,\ldots ,X_s]}{(X_1^{\ell}, \ldots ,
X_s^{\ell})} \right).
\]
\end{thm}
\begin{proof}
Let's write $Z_g$ for $\frac{ Z(\EO)}{
\mathfrak{m}_gZ(\EO )}$. The equivalence in the first sentence is  a
consequence of Theorem \ref{simplefun}(b)(ii). By Theorem \ref{maxfactor} there is an algebra isomorphism
\[
\EO (g) \cong \mathrm{Mat}_{\ell^N}\left( Z_g\right).
\]
We have isomorphisms
\begin{eqnarray*}
Z_g &\cong & \frac{Z_q\otimes_{Z_0\cap Z_q} Z_0}{\mathfrak{m}_g(Z_q\otimes_{Z_0\cap Z_q} Z_0)}  \\
&\cong & Z_q \otimes_{Z_0\cap Z_q} k_{\mathfrak{m}_g} \\
&\cong & \frac{Z_q}{(\mathfrak{m}_g\cap Z_q)Z_q}.
\end{eqnarray*}
Here $\mathfrak{m}_g\cap Z_q\triangleleft Z_0\cap Z_q$ is, in the notation of Lemma \ref{sing1}, specified by $b_{\varpi_i}^{\ell}(g)=b_i$ and $c_{\varpi_i}^{\ell}(g) = c_i$. Recalling our decomposition in Remark \ref{centre}(3),
\[
Z_q = \bigotimes_{i=1}^r Z_q^i,
\]
we see that $Z_g$ is the tensor product of rings $R_i$, for $1\leq i\leq r$, where
\begin{equation}
\label{jess}
R_i:= \frac{Z_q^i}{(\alpha_i(0)-b_i , \alpha_i(\ell) -c_i)Z_q^i}.
\end{equation}
Let's describe the possible structure of the rings $R_i$. Since $\EO (g)$ is Azumaya it follows from Proposition \ref{sing2} that we never have $b_i = 0 = c_i$. There are only two cases to consider.
\newline
\noindent
(i) $b_i \neq 0 \neq c_i$: in this case we have an algebra isomorphism
\begin{equation}
\label{eq:ss} R_i \cong \frac{\mathbb{C}[X_i]}{(X_i^{\ell}-1)}
\cong \mathbb{C}\oplus \cdots \oplus \mathbb{C}.
\end{equation}
Indeed sending $X_i \longmapsto \alpha_i(1)$ produces an isomorphism
\[
\frac{\C [X_i]}{(X_i^{\ell} - b_i^{\ell -1}c_i)} \cong R_i.
\]
Since this is a semisimple algebra of dimension $\ell$, the isomorphisms in (\ref{eq:ss}) are clear.
\newline
\noindent
(ii) $b_i \neq 0 = c_i$ (the case $b_i = 0\neq c_i$ is the same by symmetry): in this case we have an algebra isomorphism
\begin{equation}
\label{eq:trun} R_i \cong \frac{\C[X_i]}{(X_i^{\ell})}.
\end{equation}
Again sending $X_i \longmapsto \alpha_i(1)$ yields the required isomorphism.
To complete the theorem, recall that $\EO (g)$ has exactly $\ell ^{r-s}$ simple modules. But $R_i$ has exactly $\ell$ simple modules in case (i) and a unique simple module in case (ii). Therefore
\begin{eqnarray*}
\frac{Z(\EO)}{  \mathfrak{m}_gZ(\EO )} \cong \bigotimes_{i=1}^r R_i &\cong &\bigotimes_{i=1}^{s} \frac{\C[X_i]}{(X_i^{\ell})} \otimes \bigotimes_{i=s+1}^r \frac{\C[X_i]}{(X_i^{\ell} - 1)} \\
&\cong & \bigoplus_{1}^{\ell^{r-s}} \frac{\C[X_1, \ldots ,X_s]}{(X_1^{\ell},\ldots ,X_s^{\ell})}.
\end{eqnarray*}
\end{proof}
\begin{rem}
1. When $\EO (g)$ is Azumaya then, by Theorem \ref{azmor}, the
complexity of $\EO (g)$ equals $2N - \ell(w_1) - \ell(w_2)$, as
shown in \cite{gor4}.\\
2. By \cite[Theorem 2.8]{brogor1}, the maximal ideals $\mathfrak{m}_g$ of
$\mathcal{O}[G]$ which are unramified in $Z(\EO)$ form a (proper)
subset of those for which $\EO(g)$ is Azumaya, and one can read
off at once from Theorem \ref{azmor} that this set consists of
those $\mathfrak{m}_g$ with $g$ in $X_{w_0,w_0}$.
\end{rem}
\section{Representation type}
\subsection{}
\label{lessres}
To obtain general results in this section we use Theorem \ref{azmor} to remove the restriction on $\ell$ in Proposition \ref{description}.
\begin{lem}
The statement of Proposition \ref{description} is valid without the
restriction $\ell > h$.
\end{lem}
\begin{proof}
For $b\in X_{w_0,e}$ this follows from Theorem (a)(ii). For
$b\in X_{w,e}$ with $\ell (w) = N-1$ the only point in \cite{gor4} where the bound $\ell >h$ was required
was to deduce that the algebra $\UP (b)$ is Nakayama, that is its
projective indecomposable modules are uniserial. If we knew this to be
so
without the bound then the lemma would follow.
Let $b\in X_{w,e}$ and $b'\in X_{e,w_0}$ be unipotent and let
$g=b'b\in X_{w,w_0}$. By Proposition \ref{bigcell} there is an isomorphism
\[
\EO (g) \ast \mathbb{Z}_{\ell}^r \longrightarrow \UM (b') \otimes \UP
(b).
\]
As $s(w_0w)=1$ the algebra $\EO (g)$ is, by Theorem \ref{azmor}, a
truncated polynomial ring in one variable. In particular it is a
Nakayama algebra. By \cite[Theorems 1.1 and 1.3(g)]{rierie} a skew
group extension over $\mathbb{C}$ of a Nakayama algebra is again
Nakayama. Therefore $\UM(b') \otimes \UP (b)$ is a Nakayama
algebra. By definition $b'\in X_{e,w_0}$ so by the first sentence
of this proof $\UM (b')$ is a semisimple algebra, implying that
the tensor product $\UM(b')\otimes \UP(b)$ is a direct product of
matrix algebras with coefficients in $\UP (b)$. Hence $\UP (b)$ is
Morita equivalent to a Nakayama algebra and so must be a Nakayama
algebra itself. This proves the lemma.
\end{proof}
\subsection{}
\label{comptwo}
We require a general lemma from the theory of finite dimensional
algebras.
\begin{lem}
\label{repskew}
Let $S$ be a finite dimensional algebra over $\C$. Let $G$ be a finite abelian group acting by automorphisms on $S$. Then $S$ and the skew group algebra $S\ast G$ have the same representation type.
\end{lem}
\begin{proof}
Suppose we have an inclusion of algebras $S\subseteq T$. Suppose
further that $T$ has a $S$-bimodule decomposition $T=S\oplus M$. Then,
by \cite[Proposition 2]{bondro}, the representation type of $S$ is a
lower bound for the representation type of $T$ (where finite is less
than tame is less than wild).
It's clear that $S$ is a bimodule direct summand of $S\ast G$. The
character group of $G$, say $H$, acts naturally on $S\ast G$ by
\[
\chi (sg) = \chi (g) sg,
\]
and by \cite[Corollary 5.2]{rierie} the algebras $S$ and $(S\ast G)\ast H$ are Morita equivalent (this uses the fact that $|G|$ is invertible in $\C$). Combining this with the previous paragraph yields the lemma.
\end{proof}
\subsection{}\label{reptypethm}
The following lemma is the key to determining the representation
type of the algebras $\EO (g)$. Note that its validity doesn't
require that $\ell > h$.
\begin{lem}
\label{comp2} Let $w_1,w_2\in W$ and suppose
$g\in X_{w_1,w_2}$. If $\ell (w_1)+\ell (w_2) = 2N-2$ then the
algebra $\EO (g)$ is wild.
\end{lem}
\begin{proof}
The proof is based on the fact that the truncated polynomial algebra
\[
\frac{\C[X,Y]}{(X^{\ell},Y^{\ell})}
\]
has wild representation type if $\ell \geq 3$. This is a consequence of \cite[1.1(c), 1.2]{rin}.
In order that $\ell (w_1)+\ell (w_2)=2N -2 $ we must have one of the following for $1\leq i\leq r$ and $1\leq j\leq r$:
\begin{enumerate}
\item
$w_1 = w_0s_is_j , w_2 = w_0$ where $i\neq j$ (and the symmetric case obtained by exchanging the roles of $w_1$ and $w_2$);
\item
$w_1 = w_0s_i , w_2=w_0s_j$ where $i\neq j$;
\item
$w_1 = w_0s_i, w_2=w_0s_i$.
\end{enumerate}
It follows from Theorem \ref{simplefun} that Cases 1 and 2 are Azumaya so, by Theorem \ref{azmor},
\[
\EO (g) \sim_{\text{Mor}} \bigoplus^{\ell^{r-2}} \frac{\C[X, Y]}{(X^{\ell} , Y^{\ell})}.
\]
Therefore $\EO (g)$, and each of its blocks, has wild representation type.
Suppose we are in Case 3. By Proposition \ref{bigcell} we can assume without loss of generality that there is an isomorphism
\[
\gamma_g : \EO (g)\ast \mathbb{Z}_{\ell}^r \longrightarrow \UM (b')\otimes \UP (b),
\]
where
\[
g= b'b,
\]
for $b\in X_{w_0s_i,e}$ and $b'\in X_{e,w_0s_i}$ unipotent. By Lemma
\ref{repskew}, it is enough to show that $\UM (b')\otimes \UP (b)$ is
wild. By Proposition \ref{description}(b) and Lemma \ref{lessres} the algebras $\UM (b')$ and $\UP (b)$ are isomorphic to direct sums of either
\[
\text{Mat}_s \left( \frac{\C[X]}{(X^{\ell})} \right) ,
\]
or
\[
\text{Mat}_t \left( \overline{U_{\epsilon}^{\geq 0}}(\mathfrak{sl}_2)\right) .
\]
Since $\text{Mat}_s(A)\otimes \text{Mat}_t(B)\cong \text{Mat}_{st}(A\otimes B)$ it therefore suffices to show that the following algebras are wild:
\newline
(i) $\frac{\C[X,Y]}{(X^{\ell}, Y^{\ell})};$
\newline
(ii) $\overline{U_{\epsilon}^{\geq 0}}(\mathfrak{sl}_2)\otimes \frac{\C[X]}{(X^{\ell})}$;
\newline
(iii) $\overline{U_{\epsilon}^{\geq 0}}(\mathfrak{sl}_2)\otimes \overline{U_{\epsilon}^{\geq 0}}(\mathfrak{sl}_2)$.
It's clear, however, that the algebra in (ii) (respectively in
(iii)) is a skew group ring with coefficient ring
$\C[X,Y]/(X^{\ell},Y^{\ell})$ and group $\mathbb{Z}_{\ell}$
(respectively $\mathbb{Z}_{\ell}^2$). Applying Lemma \ref{repskew}
again and the comments in the first paragraph of this proof shows that
these are indeed wild.
\end{proof}
\subsection{}
\label{degen}
We need a couple of definitions from the theory of finite dimensional algebras, \cite{gab} and \cite[Chapter II]{kra}. Let
\[
\text{Bil}(n) = \{ \text{bilinear maps }m:\C^n\times \C^n \longrightarrow \C^n\}\cong \mathbb{A}^{n^3},
\]
and
\[
\text{Alg}(n) = \{\text{associative, bilinear $m$ which have an identity}\}\subseteq \text{Bil}(n).
\]
As discussed in \cite{gab} $\text{Alg}(n)$ is an affine variety,  locally closed in $\text{Bil}(n)$. The group $GL(n)$ acts on $\text{Alg}(n)$, the orbits being isomorphism classes of $n$ dimensional algebras. We let $\mathcal{O}_A$ denote the orbit in $\text{Alg}(n)$ of algebras isomorphic to $A
$. We say that \textit{$A'$ is a degeneration of $A$} if $\mathcal{O}_{A'} \subseteq \overline{\mathcal{O}_A}$, the closure of $\mathcal{O}_A$.
\begin{lem}
Let $g\in X_{w_1,w_2}$ and $g'\in \overline{X_{w_1,w_2}}$. Then
$\EO (g')$ is a degeneration of $\EO (g)$.
\end{lem}
\begin{proof}
By Proposition \ref{free} $\EO$ is a free $\mathcal{O}[G]$-module of rank $t=\ell^{\dim G}$. Let $\{x_1, \ldots, x_t\}$ be a basis for this module and define $c_{ij}^k\in \mathcal{O}[G]$ for $1\leq i,j,k\leq t$ by the following equations,
\[
x_ix_j = \sum_k c_{ij}^k x_k.
\]
Then for any $g\in G$ the structure constants of $\EO (g)$ with respect to the basis $\{ x_ i + \mathfrak{m}_g\EO \}$ are given by $(g(c_{ij}^k)) = (c_{ij}^k(g))$. As a result the map
\[
\alpha : G\longrightarrow Alg(t) \subseteq \mathbb{A}^{t^3},
\]
defined by $\alpha (g) = ( c_{ij}^k (g))$, is a morphism of varieties.
Let $g\in X_{w_1,w_2}$ and $g'\in \overline{X_{w_1,w_2}}$. By Theorem \ref{simplefun} $X_{w_1,w_2}$ is a dense open set of $\overline{X_{w_1,w_2}}$ over which all algebras in the family $(\alpha(z))_{z\in \overline{X_{w_1,w_2}}}$ are isomorphic to $\EO (g)$. It follows from \cite[Proposition 3.5
and Section 3.7]{kra} that $A_{g'}\cong \EO (g')$ is a degeneration of $\EO (g)$.
\end{proof}
\begin{rem}
The above proof is also valid for the reduced quantum Borels.
Namely, if $b\in X_{w,e}$ and $b'\in \overline{X_{w,e}}$ then $\UP
(b')$ is a degeneration of $\UP (b)$.
\end{rem}
\subsection{}
\label{reptypethm2}
The following statement was proved by Gei{\ss} in \cite{gei},
\begin{equation}
\label{geiss} \text{ Let $A'$ be a degeneration of $A$. If $A$ is
wild then so is $A'$.}
\end{equation}
This allows us to complete the classification of the
representation type of the algebras $\EO (g)$, without the
restriction $\ell > h$.
\begin{thm}
Let $\ell$ be good. Let $w_1, w_2 \in W$ and suppose that $g\in X_{w_1,w_2}$.
(i) If $\ell(w_1) + \ell(w_2)\geq 2N - 1$ then $\EO (g)$ has finite representation type.
(ii) If $\ell(w_1) + \ell(w_2) \leq 2N - 2$ then $\EO (g)$ has wild representation type.
\end{thm}
\begin{proof}The first part follows directly from Theorem \ref{allaz}. Indeed, if
$\ell(w_1)+\ell(w_2) = 2N$ then the algebra $\EO (g)$ is semisimple whilst if
$\ell(w_1)+\ell(w_2)=2N-1$ then $\EO (g)$ is Morita equivalent to a direct sum of truncated
 polynomial algebras in one variable, hence Nakayama.
Let $\preccurlyeq$ denote the Bruhat-Chevalley order on $W$. For
the second part note that if $\ell(w_1) + \ell(w_2)< 2N - 2$ then
there exists $u, v\in W$ such that $w_1\preccurlyeq u$,
$w_2\preccurlyeq v$ and $\ell(u) + \ell(v) = 2N -2$. Arguing
exactly as in \cite[Theorem 2.1]{rich2} it follows that
\[
\overline{X_{u,v}} = \overline{BuB}\cap \overline{B^-vB^-} =
\coprod_{u'\preccurlyeq u,v'\preccurlyeq v} X_{u',v'},
\]
so $g\in \overline{X_{u,v}}$. Let $g'\in X_{u,v}$. Then by Lemma \ref{reptypethm} $\EO (g')$ is wild and by Lemma \ref{degen} $\EO (g) $ is a degeneration of $\EO (g')$. Therefore, by (\ref{geiss}), $\EO(g)$ must be wild.
\end{proof}
\subsection{}
\label{utype}
We now tackle representation type for $\UP (b)$.
\begin{cor}
Let $\ell$ be good. Let $w\in W$ and suppose that $b\in X_{w,e}$.
(i) If $\ell (w) \geq N-1$ then $\UP (b)$ has finite representation type.
(ii) If $\ell(w) \leq N-2$ then $\UP (b)$ has wild representation type.
\end{cor}
\begin{proof}The first part follows from Lemma \ref{lessres} and the observation that $\overline{U_{\ep}^{\geq 0}(\mathfrak{sl}_2)}$ has finite representation type by Lemma \ref{repskew}. For (ii), arguing as in the proof of Theorem \ref{reptypethm2} we see that it is sufficient to show that $\UP (b)$ is wild in the case that $\ell(w) = N-2$, that is $w = w_0s_is_j$ for $i\neq j$. Let $g\in X_{w_0s_is_j,w_0}$ be such that $g=b'b$
 where $b'\in X_{e,w_0}$ and $b\in X_{w_0s_is_j,e}$ are unipotent. By
Proposition \ref{bigcell} we have an algebra isomorphism
\[
\EO (g)\ast \mathbb{Z}_{\ell}^r \cong \UM (b') \otimes \UP (b).
\]
Since $b'\in X_{e,w_0}$ we have, by Proposition \ref{description}(a)
and Lemma \ref{lessres}, an isomorphism
\[
\UM (b') \cong
\bigoplus^{\ell^{r-s(w_0)}}\mathrm{Mat}_{\ell^{\frac{1}{2}(N+s(w_0))}}(\C).
\]
Therefore
\[
\EO (g)\ast \mathbb{Z}_{\ell}^r \cong
\bigoplus^{\ell^{r-s(w_0)}}\mathrm{Mat}_{\ell^{\frac{1}{2}(N+s(w_0))}}\left(
\UP (b) \right).
\]
So $\UP (b)$ has the same representation type as $\EO (g)\ast
\mathbb{Z}_{\ell}^r$ and hence, by Lemma \ref{repskew}, as $\EO
(g)$. Now apply Theorem \ref{reptypethm2}.
\end{proof}
\begin{rem}
When $b$ is the identity element of $B^-$ the results of the corollary were obtained in \cite{cib}.
\end{rem}
\section{Algebras with group actions}
\subsection{}
\label{skewgroup}
Let $R$ be a finite dimensional $k$-algebra and $G$ be a finite
abelian group. Assume that the
characteristic of $k$ is prime to the order of $G$. Suppose $G$ acts
as algebra automorphisms on $R$, so we have a group homomorphism
\[
\tau :G \longrightarrow \text{Aut}_{\text{k-alg}}(R).
\]
If $M$ is a finite dimensional $R$-module we let $^g\!M$ denote
the $R$-module whose underlying abelian group is $M$ and whose
action is given by $r\cdot m = \tau(g)^{-1}(r)m$.  Given $g\in G$
there is a functor
\[
F_g : R-mod \longrightarrow R-mod,
\]
which takes $M$ to $^g\!M$ and sends $f:M\longrightarrow N$ to
$f:\: ^g\!M\longrightarrow \: ^g\!N$.
Fix a simple $R$-module $V$ and let $V(g) = \: ^g\!V$. Throughout
this section we shall assume that
\[
\{ V(g): g\in G \} \textit{ is a complete set of non-isomorphic
simple $R$-modules.}
\]
In particular, this assumption implies that $\tau$ is a
monomorphism, and that the simple modules share a fixed
$k$-dimension, $t$ say. Fix $P$, a projective cover of $S$, and
let $P(g) = \: ^g\!P$. By the above assumption
\[
Q \cong \bigoplus_{g\in G}P(g)
\]
is a projective generator for $R-mod$. Let $E=\text{End}_R(Q)$.
Given $g\in G$ let $\sigma_g :G \longrightarrow G$ denote the left
regular action, that is $\sigma_g(h)=gh$. Considering elements of
$Q$ as ordered $|G|$-tuples of elements of $P$ we can define
$\psi_g
:Q \longrightarrow Q$ to be the additive map which acts as the
permutation $\sigma_g$ on the $|G|$-tuple. In other words an
element concentrated in the $h^{\text{th}}$ position is sent to
the $gh^{\text{th}}$ position.
\begin{lem}
For $g\in G$ let $\psi_g:Q \longrightarrow Q$ be as above. Then
\\
(i) $\psi_g(r\cdot q) = \tau(g)(r)\cdot \psi_g(q)$ for all $r \in
R$ and $q \in Q$;
\\
(ii) $\psi_g\psi_h = \psi_{gh}$.
\end{lem}
\begin{proof}
The second claim is obvious. For the first we can assume that $q$ is
concentrated in the $h^{\text{th}}$ position. Then
\[
\psi_g(r\cdot q) = \psi_g(\tau(h)^{-1}(r)q) =
\tau(h)^{-1}(r)q,
\]
where the right hand side is concentrated in the $gh^{\text{th}}$
position. On the other hand, since $\psi_g(q)$ is non-zero only in
the $gh^{\text{th}}$ position, we have
\[
\tau(g)(r)\cdot \psi_g(q) = \tau(gh)^{-1}(\tau(g)(r))q =
\tau(h)^{-1}(r)q,
\]
as required.
\end{proof}
\subsection{}
\label{skewgr}
 For $g\in G$ let $\tilde{\tau}(g):E\longrightarrow
E$ send $\phi$ to the map $\psi_g\circ \phi\circ \psi_g^{-1}$. The
lemma ensures that this is a well-defined $k$-algebra automorphism
and that the induced map
\[
\tilde{\tau}: G \longrightarrow
\text{Aut}_{\text{$k$-alg}}(E^{\text{op}})
\]
is a group homomorphism. Repeating the comments of the first paragraph
of this section we have a functor
\[
\tilde{F}_g : E^{\text{op}}-mod \longrightarrow E^{\text{op}}-mod,
\]
sending $N$ to $^gN$ and fixing homomorphisms.
Observe that $Q$ is an $(R, E^{op})$-bimodule with $r\cdot q \cdot \phi
= \phi(r\cdot q)$ for all $r\in R$, $q\in Q$ and $\phi \in E$. There is an equivalence of categories
\[
R-mod \longrightarrow E^{op}-mod
\]
given on objects by sending $M$ to $\hm{R}{Q}{M}$. The inverse
equivalence sends $N$ to $Q\otimes_{E^{\text{op}}}N$. This
equivalence induces two functors for each $g\in G$, namely
\[
\alpha_g, \beta_g : R-mod \longrightarrow R-mod
\]
where $\alpha_g(Q\otimes_{E^{\text{op}}}N) =
^g\!Q\otimes_{E^{\text{op}}}N$ and
$\beta_g(Q\otimes_{E^{\text{op}}}N) =
Q\otimes_{E^{\text{op}}}{\!^gN}$. So $\alpha_g$ corresponds to $F_g$
and $\beta_g$ to $\tilde{F}_g$.
\begin{prop}
The functors $\alpha_g$ and $\beta_g$ are naturally isomorphic.
\end{prop}
\begin{proof}
Let $\theta_g : Q\otimes_{E^{\text{op}}}N \longrightarrow
^{g^{-1}}\!Q\otimes_{E^{\text{op}}} \!^gN$ send $q\otimes n$ to
$\psi_g(q)\otimes n$. We must check this is well-defined. First
note that, for $\phi \in E^{\text{op}}, \psi_g(q \phi) = \psi_g(q)
\tilde{\tau}(g)(\phi)$. Thus we have, for $q \in Q$ and $n \in N$,
\begin{eqnarray*}
\theta_g(q \phi \otimes n)& =& \psi_g(q\phi ) \otimes n \\ & = &
\psi_g(q)\tilde{\tau}(g)(\phi)\otimes n \\ & = & \psi_g(q)\otimes
\tilde{\tau}(g)(\phi)\cdot n \\ & = & \psi_g(q)\otimes \phi n \\ &
=& \theta_g(q\otimes \phi n) .
\end{eqnarray*}
Moreover $\theta_g$ is an $R$-module isomorphism, since, for $r
\in R$,
\begin{eqnarray*}
\theta_g(rq\otimes n)& =& \psi_g(rq)\otimes n \\
&=&
\tau(g)(r)\psi_g(q)\otimes n\\
 & =& r\cdot \psi_g(q) \otimes n\\
& = & r\cdot \theta_g(q\otimes n).
\end{eqnarray*}
Since $\theta_g$ is natural in $N$ it follows that
$\alpha_g^{-1}\beta_g$ is naturally isomorphic to the
identity functor. One shows similarly that $\beta_g^{-1}\alpha_g $ is also
naturally isomorphic to the identity functor.
\end{proof}
\indent Given $g\in G$ let $\pi_g \in E$ be the primitive
idempotent corresponding to projection onto $P(g)$ followed by the
canonical injection of $P(g)$ into $Q$. For $h\in G$ it is easy to
see that we have $\tilde{\tau}(h)(\pi_g) = \pi_{hg}$. Since
${}_{E^{\text{op}}}^{g}\!Q \cong {}_{E^{\text{op}}}Q$ for all $g \in G$,
one sees that $Q$ is a free $E^{\text{op}}$-module of rank $t$ for
some $t \geq 1$, so that
\begin{eqnarray*}
& \textit{$R \cong \mathrm{Mat}_{t}(S)$ where $S$ is a basic algebra
on which $G$ acts permuting}\\
& \textit{a set of minimal primitive idempotents
simply transitively.}
\end{eqnarray*}
Notice that $t$ as it appears in the above statement coincides with its earlier definition as the (shared) dimension of the simple $R$-modules. We let $\{ e_g: g\in G\}$ be the above set of minimal primitive
idempotents of $S$, and let $X=X(G)$ be the character group of
$G$. Since $G$ is abelian we have a decomposition
\[
S = \bigoplus_{\chi\in X} S_{\chi}
\]
where $S_{\chi}= \{ s \in S : \tau(g)(s) = \chi(g)s\}$, so that
$S$ is an $X$-graded algebra. Given $\chi \in X$ we define
\[
y_{\chi} = \sum_{g\in G} \chi^{-1}(g)e_g\in S_{\chi}.
\]
If the exponent of $G$ is $\ell$ then we find
\[
y_{\chi}^{\ell} = \sum_{g\in G}\chi^{-1}(g^{\ell})e_g^{\ell} =
\sum_{g\in G}e_g = 1;
\]
moreover, for $\chi, \eta \in X$ with $\chi \neq \eta$,
\[
y_{\chi}y_{\eta} = y_{\chi\eta}.
\]
Thus $y_{\chi}$ is a unit in $S_{\chi}$ for $\chi \in X$, and
$\sum_{\chi \in X}k y_{\chi}$ is a subalgebra of $S$
normalising $S_1$, and isomorphic to $k X$ and hence to $k G$.
\begin{thm}
Retain the notation and hypotheses of the above paragraphs. There is an isomorphism $S\cong S_1 \ast G$, where the right hand side is a skew group ring. Moreover $S_1$ is
scalar local. Thus $R \cong \mathrm{Mat}_{t}(S_{1} \ast G) \cong
\mathrm{Mat}_{t}(S_1) \ast G$.
\end{thm}
\begin{proof}
The discussion above shows that $R$ is a skew group ring $R_{1}\ast
G$, and that we may reduce to the case where $R$ is a basic
algebra. By \cite[Theorem 4.2]{Pas} $J(R_1)R = J(R)$. By Lying Over for $R_1
\subseteq R$, \cite[Theorem 16.6]{Pas}, we deduce that $R_1$ is a basic algebra. Thus
$R_1/J(R_1)$ is a finite direct sum of copies of $k$. Commutativity of $R/J(R)$ forces the action of $G$ on $R_{1}/J(R_{1})$ to be
trivial. Since there are exactly $|G|$ simple $R$-modules it follows
that $R_1$ is (scalar) local.
\end{proof}
\subsection{Blocks and quivers}
\label{blquiv} The blocks and quiver of $S$, (and hence of $R$),
are determined by the conjugation action of its subgroup $G$ on
$J(S_1)/J(S_1)^2$. We make this statement precise through the
Morita equivalence of Paragraph \ref{skewgr} and the following
result. See \cite{ARS} for the terminology used here, recalling
also the definition of a multiply-edged Cayley graph from
\ref{I3}.
\begin{prop}
Let $T$ be the skew group ring $T_1 \ast G$ of a finite abelian
group $G$ over the scalar local finite dimensional $k$-algebra
$T_1$, with $k$ algebraically closed of characteristic not
dividing $|G|$. Let $J$ be the Jacobson radical of $T_1$, so
$J/J^2$ is a $kG$-module under the conjugation action
\begin{eqnarray*}
g.(t + J^2) \quad = \quad gtg^{-1} + J^2,
\end{eqnarray*}
for $g \in G$ and $t \in J$. Let
\begin{eqnarray}
\label{decomp}
J/J^2 \quad = \quad \sum_{\chi \in
X(G)}^{\oplus}V_{\chi}^{(m_{\chi})}
\end{eqnarray}
be the decomposition of $J/J^2$ as a direct sum of irreducible
$kG$-modules under this action. Define
\begin{eqnarray*}
Y \quad := \quad Y(G,T) \quad = \quad <\chi : m_{\chi} \neq 0> \quad \subseteq \quad X(G),
\end{eqnarray*}
and
\begin{eqnarray*}
D \quad := \quad C(G,T) \quad = \quad C_G (J/J^2) \quad = \quad \{
g \in G : \chi(g) = 1 \textit{ for all } \chi \in Y \}.
\end{eqnarray*}
\indent (i) The quiver $Q_T$ of $T$ has vertices $\{v_{\chi} :
\chi \in X(G) \}$, and an arrow $v_{\eta} \longrightarrow v_{\mu}$
if and only if $\mu = \eta \chi^{-1}$ for some $\chi$ with
$m_{\chi} \neq 0$.\\ \indent (ii) The number of blocks of $T$ is
$|X(G) : Y| = |D|$.\\ \indent (iii) Each block of $T$ has an
identical quiver, namely the multiply-edged Cayley graph of $Y$
with respect to the generating set $\{ \chi^{-1} : m_{\chi} \neq 0
\}$ of $Y$, with $m_{\chi}$ copies of the edge $\chi^{-1}$
starting at each vertex .\\ \indent (iv) The following statements
are equivalent:\\ \indent (a) The blocks of $T$ are trivial.\\
\indent (b) $Y = \{1\}$.\\ \indent (c) $G$ centralises $J/J^2$.\\
\indent (d) $T$ is the ordinary group ring $T_1G$.
\end{prop}
\begin{proof}
(i) Since $J(T) = J \ast G$ by \cite[Theorem 4.2]{Pas}, $T/J(T)
\cong kG$, so that $T$ is basic and has quiver with vertices
labelled by $X(G)$. By definition (see e.g.\cite[page 65]{ARS}),
to find the arrows of the quiver of $T$ we may assume without loss
that $J^2 = 0$. The orthogonality relations for $kG$ show that the
primitive idempotents of $kG$ (and hence of $T$) are $\{e_{\chi} =
1/|G| \sum_{g \in G} \chi(g^{-1})g : \chi \in X(G) \}$. Taking a
basis for $J$ consisting of eigenvectors $v_{\chi} \in V_{\chi}$
with respect to its structure (\ref{decomp}) as $kG$-module, and
noting that $e_{\mu}v_{\chi}e_{\eta}$ is non-zero if and only if
$\mu = \eta \chi^{-1}$, one finds $e_{\mu}Je_{\eta}$ is non-zero
if and only if $\mu = \eta \chi^{-1}$, proving (i).\\ \indent (ii)
It's clear from (i) that two vertices $v_{\eta}$ and $v_{\mu}$ are
in the same connected component of the quiver if and only if
$\eta$ and $\mu$ are in the same coset of $Y$ in $X(G)$. The final
equality is obvious, since $Y = X(G/D)$.\\ \indent (iii) Immediate
from (i).\\ \indent (iv) (a) $\Longrightarrow$ (b): By (ii).\\
\indent (b) $\Longrightarrow$ (c): By definition of $Y$.\\ \indent
(c) $\Longrightarrow$ (d): It's easy to show that any choice of
lifts to $J$ of a basis of $J/J^2$ generate $J$ as a $T_1$-module
- see for example \cite[Theorem III.1.9(a)]{ARS}. Thus, if (c)
holds then $G$ operates unipotently on $J$ and hence on $T_1$. The
assumption on the characteristic of $k$ ensures that the action of
$G$ on $T_1$ is completely reducible. So (d) follows.\\ \indent
(d) $\Longrightarrow$ (a): Trivial.
\end{proof}
\begin{rem}
For the applications of Proposition \ref{blquiv} below it's
convenient to formulate the following easy generalisation. Namely,
suppose that $T$ is a finite dimensional $k$-algebra containing a
group algebra $kG$ of a group $G$ whose order is invertible in
$k$, such that the primitive central idempotents of $T$ are
precisely the primitive idempotents of $kG$. Then conclusions (i),
(ii) and (iii) of the proposition remain true, with $J(T)$
replacing $J$.
\end{rem}
\section{Reduced quantum Borels}
\subsection{}
\label{qbors}
Let $w \in W$ and choose $b \in X_{w,e}$. From Lemma
\ref{ellform}(ii) and the theory developed in Section 5, we know
that $\UP (b)$ is a matrix ring over a skew group ring whose
coefficient ring is scalar local. We proceed now to identify the
components of this structure.\\
\indent Let $w = s_{i_1}s_{i_2} \ldots s_{i_t}$ be a reduced expression
for $w$ as a product of simple reflections, so
\begin{eqnarray}
\label{wlength}
\ell (w) = t,
\end{eqnarray}
and let $\beta_{i_1}, \dots , \beta_{i_t}$ be the corresponding
ordered subset of the positive roots of $\g$, with corresponding
PBW-type generators $E_{\beta_{i_1}}, \dots , E_{\beta_{i_t}}$ in
$\UP$, \cite[Theorem 8.24]{janqg}, which we renumber respectively
as $\beta_1, \ldots , \beta_t$ and as $E_{\beta_1}, \ldots ,
E_{\beta_t}$. Writing $w_0 = ww_1$ for an element $w_1$ of $W$
with $\ell (w) + \ell (w_1) = N = \ell (w_0)$, we obtain
corresponding PBW-type elements $E_{\beta_1}, \ldots ,
E_{\beta_t}, E_{\beta_{t+1}}, \ldots , E_{\beta_N}$ of
$U^{>0}_{\epsilon}$. Set $A(b)$ to be the subalgebra of
$U^{>0}_{\epsilon} (b)$ generated by the images in the latter
algebra of $E_{\beta_1}, \ldots , E_{\beta_t}$. (We'll abuse
notation by using the same notation for the image of $E_{\beta_i}$
in $A(b)$, for $i = 1, \dots , t$, and also in
$U^{>0}_{\epsilon}(b)$, for $i = 1, \ldots , N$.)
Now
\begin{eqnarray*}
\UP (b)  = U^{> 0}_{\epsilon}(b)\ast P_{\ell},
\end{eqnarray*}
a skew group ring. By the PBW-type theorem,
\begin{eqnarray*}
\mathrm{dim}_{\C}(U^{> 0}_{\epsilon}(b))  =  \ell^{N},
\end{eqnarray*}
with basis $\{E_{\beta_1}^{m_1}E_{\beta_2}^{m_2} \ldots
E_{\beta_N}^{m_N} : 0 \leq m_i < \ell \}$. In a similar way,
$
\mathrm{dim}_{\C}(A(b)) = \ell^t.
$
Note that $A(b)$ is normalised by $P_{\ell}$ in its conjugation
action on $U^{> 0}_{\epsilon}(b)$, so that $A(b)\ast P_{\ell}$ is
a skew group subalgebra of $\UP (b)$, and
\begin{eqnarray}
\label{one} \mathrm{dim}_{\C}(A(b)\ast P_{\ell}) = \ell^{t + r}.
\end{eqnarray}
By \cite[Theorems 4.4 and 4.5]{decpro49} and \cite[Corollary
2.8]{gor4} $A(b)$ and $A(b)\ast P_{\ell}$ are
semisimple, so that
\begin{eqnarray}
\label{two}
 A(b)\ast P_{\ell} \quad \cap \quad J(\UP (b)) \quad =  \quad 0.
\end{eqnarray}
\indent By Theorem \ref{appnot}(i),
\begin{eqnarray}
\label{three} \mathrm{dim}_{\C}(\UP(b)/J(\UP (b))\quad = \quad \ell^{r -
s(w)}(\ell^{1/2(t + s(w))})^2 \quad = \quad \ell^{r+t}.
\end{eqnarray}
We conclude from (\ref{one}), (\ref{two}) and (\ref{three}) that
\begin{eqnarray}
\label{four}
 \UP(b)/J(\UP (b))\quad \cong \quad A(b)\ast P_{\ell}.
\end{eqnarray}
Notice that (\ref{four}) explains the correspondence between
simple $\UP (b)-$ and $A(b)\ast P_{\ell}-$modules, \cite{decpro49}.
\subsection{}
\label{skewgprgs} Thanks to Lemma \ref{ellform}(i) there is a
canonical isomorphism between $Q_{\ell}$ and the character group
of $P_{\ell}$, given by $\alpha \longmapsto
\epsilon^{(\alpha,-)}$. Thus $P_{\ell}$ and $Q_{\ell}$ are isomorphic
groups, but we shall nevertheless use the two notations to denote
(respectively) the group of automorphisms of $\UP(b)$ induced by
conjugation by its subgroup $P_{\ell}$ of units, and the action of
$Q_{\ell}$ by right winding automorphisms.
\begin{lem}
\label{Abideals}
 Keep the notation as above and in the previous
subsection.\\ (i) $A(b)$ acts faithfully on each simple
$\UP(b)$-module.\\ (ii) The only ideals of $A(b)$ which are
$P_{\ell}$-invariant are
 $0$ and $A(b)$.
\end{lem}
\begin{proof}
 (i) By Theorem \ref{appnot}(i) the winding automorphisms afforded by $Q_{\ell}$ permute the irreducible $U^{\geq 0}_{\epsilon}(b)$-modules transitively.
 Thanks to the definition of the coproduct and of the winding automorphisms (see (\ref{appnot})),
 the latter act  trivially on $A(b)$. Hence, if a simple $A(b)$-module does not feature in a component of any one simple $U^{\geq 0}_{\epsilon}(b)$-module,
 then that simple $A(b)$-module doesn't feature in any of the simple  $U^{\geq 0}_{\epsilon}(b)$-modules. But this would imply that
 $J(U^{\geq 0}_{\epsilon}(b)) \cap A(b) \neq 0$, contradicting (\ref{two}). So each simple $U^{\geq 0}_{\epsilon}(b)$
 -module is faithful for $A(b)$.  \\
 \indent (ii) Suppose $I$ were a non-zero proper $P_{\ell}$-invariant ideal of $A(b)$. Then, thanks
to (\ref{four}), $\hat{I} := I \ast P_{\ell} + J(U^{\geq
0}_{\epsilon}(b))$ would be a proper (two-sided) ideal of $U^{\geq
0}_{\epsilon}(b)$ , and $\hat{I}$ would then annihilate a simple
$U^{\geq 0}_{\epsilon}(b)$-module $V$. Thus $IV = 0$ contradicting
(i).
\end{proof}
\subsection{}
The following lemma is a special case of \cite[Lemmas 6.1.5,
6.1.6] {pas2}.
\begin{lem}
\label{idemperm} Let $R$ be a ring and let $1 = e_1 + e_2 + \ldots
+ e_n$ be a sum of orthogonal idempotents. Let $G$ be a subgroup
of the group of units of $R$ and assume that $G$ permutes $\{e_1,
\dots , e_n\}$ transitively by conjugation. Then $R \cong
\mathrm{Mat}_n (S)$, where $S \cong e_1 R e_1$.
\end{lem}
\subsection{}
Now let
\[
1 = e_1 + e_2 + \ldots + e_n
\]
be a decomposition of $1 \in A(b)$ as a sum of primitive central
idempotents of $A(b)$, and let
\[
H := C_{P_{\ell}}(e_1) = \{x \in P_{\ell} : xe_1 = e_1x \}.
\]
Since $P_{\ell}$ is abelian and since, by Lemma
\ref{Abideals}(ii),
\begin{eqnarray}
\label{idemtrans} P_{\ell} \textit{ acts transitively on } \{e_1,
\dots , e_n\},
\end{eqnarray}
$H  = \cap_{i=1}^{n}C_{P_{\ell}}(e_i)$. Let the dimension of a
simple $A(b)$-module be $u$, so
\[
A(b) \quad = \quad \left( \mathrm{Mat}_{u}(\C)\right) ^{\oplus (n)}.
\]
Thus $e_1( A(b)\ast P_{\ell}) e_1 \cong \mathrm{Mat}_u (\C) \ast H$,
and it follows from Lemma \ref{idemperm} that
\begin{eqnarray}
\label{2mat} A(b)\ast P_{\ell} \quad \cong \quad \mathrm{Mat}_n (
\mathrm{Mat}_u (\C) \ast H ).
\end{eqnarray}
By the Skolem-Noether theorem the action of $H$ on $\mathrm{Mat}_u (\C)$ is
by inner automorphisms. So by \cite[Proposition 12.4]{Pas},
(\ref{2mat}) yields
\begin{eqnarray}
\label{twtens}
 A(b)\ast P_{\ell} \quad = \quad \mathrm{Mat}_{nu}(\C)
 \otimes_{\C} \C^t \tilde{H},
\end{eqnarray}
for some twisted group algebra $\C^t \tilde{H}$ of a group
$\tilde{H}$ isomorphic to $H$. Thus $t :\tilde{ H} \times \tilde{H} \longrightarrow
\C^{\ast}$
is a $2$-cocycle, for which we set
\[
\tilde{H}_0 \quad := \quad \{ h \in \tilde{H} : t(h,a) = t(a,h)
\textit{ for all } a \in \tilde{H} \},
\]
 a subgroup of $\tilde{H}$.
\begin{lem}
\label{twist}
 $Z(\mathrm{Mat}_{nu}(\C) \otimes \C^t \tilde{H} )
= \C^t \tilde{H}_0$.
\end{lem}
\begin{proof}
See \cite[Ex.3, p.176]{Pas}.
\end{proof}
\subsection{}
\label{findH0} By (\ref{twtens}), Lemma \ref{twist} and the known
parameters for $A(b)\ast P_{\ell}$ (see Theorem (\ref{appnot})(i))
,
\begin{eqnarray}
\label{H0order} |\tilde{H}_0| \quad = \quad \ell^{r-s(w)}.
\end{eqnarray}
Tracing through the isomorphisms (\ref{2mat}), (\ref{twtens}) and Lemma
\ref{twist}, one sees that the centre of $A(b) \ast P_{\ell}$,
$\C^t \tilde{H}_0$, is a group algebra of a group $\tilde{H}_0$,
where
\begin{eqnarray}
\label{tildeH}
\tilde{H}_0 \quad = \quad \{ \alpha_x x : x \in H_0 \}.
\end{eqnarray}
Here, $\C^t\tilde{H}_0$ is actually an ordinary untwisted group algebra of $\tilde{H}_0$, because it is commutative, by \cite[Lemma 1.2.9]{pas2},
and $H_0$ is a subgroup of $H$ (which itself is a subgroup of
$P_{\ell}$); and, for $x \in H_0, \alpha_x \in A(b)$ is a unit,
conjugation by which coincides with the action of $x^{-1}$ by
conjugation on $A(b)$. Our aim now is to identify $H_0$. To do so we must
\begin{eqnarray}
\label{ellhyp}
\textit{assume that $\ell$ is prime to the order of $w$.}
\end{eqnarray}
Let
\begin{eqnarray*}
X := \{ x \in Q_{\ell} : \tau_x \textit{ fixes the simple } A(b)
\ast P_{\ell}-\textit{modules}\},
\end{eqnarray*}
in the notation of \ref{appnot}. Thus $X$ fixes the primitive central
idempotents of $A(b)\ast P_{\ell}$. Note that $Q_{\ell}$, (and hence $X$), act trivially on $A(b)$, by the definition of the coproduct. Thus, in view of (\ref{twtens}), Lemma \ref{twist} and (\ref{tildeH}),
\begin{eqnarray}
\label{Xeq}
X \quad = \quad C_{Q_{\ell}}(H_0).
\end{eqnarray}
On the other hand Theorem \ref{appnot}(i) states that
\begin{eqnarray}
\label{Xiden}
X = C_{Q_{\ell}}(P_{\ell}^w).
\end{eqnarray}
Thus, from (\ref{Xeq}), (\ref{Xiden}) and Lemma \ref{ellform}(i)
we conclude that
\begin{eqnarray}
\label{X}
H_0 \quad = \quad C_{P_{\ell}}(X) \quad = \quad P_{\ell}^{w}.
\end{eqnarray}
We summarise what we have proved in the following
theorem.
\begin{thm}
\label{group} Let $b \in X_{w,e}$ and suppose (in addition to the
standing hypotheses \ref{not} on $\ell$) that $\ell$ is prime to
the order of $w$. Then the primitive central idempotents of
$\UP(b)/J(\UP(b))$ are the images of the primitive idempotents of
the subalgebra $\C \tilde{P}_{\ell}^w$ of $A(b) \ast P_{\ell}
\subseteq \UP (b)$. For $x \in P_{\ell}^w$ there is a unit
$\alpha_x$ in $A(b)$ such that $\tilde{P}_{\ell}^w = \{
\alpha_{x}x : x \in P_{\ell}^w \}$.
\end{thm}
\subsection{}
\label{booth}
We require a lemma concerning roots.
\begin{lem}
Let $w=s_{i_1}\ldots s_{i_{t}} \in W$ with $\ell(w) = t$, and let
$A_w = \sum_{j=1}^{t}\mathbb{Z}\beta_j$, in the notation of
\ref{qbors}. Then
\[
A_w = \sum_{j=1}^{t} \mathbb{Z}\alpha_{i_j}.
\]
\end{lem}
\begin{proof}
It is clear that $A_w\subseteq
\sum_{j=1}^{\ell(w)}\mathbb{Z}\alpha_{i_j}$ since by construction
$\beta_j$ is a combination of the simple roots $\alpha_{i_k}$ for
$1\leq k\leq \ell(w)$.
We prove the opposite inclusion by induction on $\ell(w) =t$, the
case $\ell (w) \leq 1$ being trivial. Let $w' = s_{i_2}\ldots
s_{i_t}$ so that $\ell (w') = t-1$, and let
$\tilde{\beta}_1,\ldots \tilde{\beta}_{t-1}$ be the set of
positive roots corresponding to $w'$. We have $\beta_1 =
\alpha_{i_1}$ and $\beta_k = s_{i_1}(\tilde{\beta}_{k-1}) =
\tilde{\beta}_{k-1} + n_k\alpha_{i_1}$ for $2\leq k\leq t$ and
$n_k\in \mathbb{Z}$. Therefore
\[
A_{w'} = \sum_{k=2}^{t}\mathbb{Z}\tilde{\beta}_{k-1} \subseteq A_w,
\]
and so $A_{w'} + \mathbb{Z}\alpha_{i_1} \subseteq A_w$. By induction
$A_{w'} = \sum_{k=2}^t\mathbb{Z}\alpha_{i_k}$ so we have $\sum_{k=1}^t
\mathbb{Z}\alpha_{i_k} \subseteq A_w$ as required.
\end{proof}
\subsection{Blocks and quiver for $\UP(b)$}
\label{Upquiv} Retain hypothesis (\ref{ellhyp}). From Theorem
\ref{group} we see that the hypotheses of Remark \ref{blquiv}
apply. To state the consequences for $\UP (b)$ it remains only to
identify the groups $G$ and $D$, or equivalently $X(G)$ and $Y$,
in Proposition \ref{blquiv}. We already know from Theorem
\ref{group} that $G = \tilde{P}_{\ell}^w$. Since
\begin{eqnarray*}
J(\UP (b)) \quad = \quad J(U_{\epsilon}^{>0}(b))\ast P_{\ell},
\end{eqnarray*}
we can write
\begin{eqnarray*}
J_1 \quad = \quad J(U_{\epsilon}^{>0}(b));
\end{eqnarray*}
then
\begin{eqnarray*}
D \quad = \quad C_{\tilde{P}_{\ell}^w}(J_1) \quad = \quad
C_{\tilde{P}_{\ell}^w}(J_1/J_{1}^2)),
\end{eqnarray*}
and we'll denote this group by $C_{\ell}^w$.
\begin{thm}
Continue with the notation and hypotheses of this section. In
particular, $b \in X_{w,e}$ with $\ell$ prime to the order of $w$,
and $w$ as in
\ref{qbors}. Let $d$ be the number of simple
reflections not occurring in a reduced expression for $w$.\\
\indent (i) The number of blocks of $\UP (b)$ is $|C_{\ell}^w|$.\\
\indent (ii) Set $B^w = (A_w)^{\perp}=\left(\sum_{j=1}^{t}
\mathbb{Z}\beta_j \right)^{\perp}$, a subgroup of $P$. Then $B^w
\subseteq P^w$ and $B_{\ell}^w
:= B^w/{\ell}B^w \subseteq \tilde{P}_{\ell}^w$, with $B_{\ell}^w \cap
C_{\ell}^w = 0$, and $|B_{\ell}^w | = \ell^d$.\\ \indent (iii) The
quiver of each block of $\UP (b)$ is a multiply-edged Cayley graph
of $\tilde{P}_{\ell}^w/C_{\ell}^w$ with respect to the generating
set given by the inverses of the weights of $\tilde{P}_{\ell}^w$
on $J_1/J_{1}^2$. In particular, $B_{\ell}^w$ embeds in (the set
of vertices of) each block.
\end{thm}
\begin{proof} Only (ii) and the final sentence of (iii) are not
immediate from Remark \ref{blquiv}. By Lemma \ref{booth},
\begin{eqnarray}
\label{bwbasis}
 B^w =
\sum \mathbb{Z}\varpi_j
\end{eqnarray}
where the sum is taken over the set $\mathcal{C}_{w}$ of all
$\varpi_j$ such that $s_j$ does not appear in a reduced expression
of $w$. Thus each such $\varpi_j$ is fixed by $w$ and we have $B^w
\subseteq P^w$. The same argument shows that $B^w_{\ell}\subseteq
P^w_{\ell}$. That $|B_{\ell}^w | = \ell^d$ is clear from
(\ref{bwbasis}).
Finally, note  that for $x \in B_{\ell}^w$, the
corresponding unit $\alpha_{x}$ of $A(b)$ is just a scalar, since
$B_{\ell}^w$ commutes with $A(b)$. Since these scalars have
trivial action on $J_1/J_1 ^2$, the final sentence of (iii)
follows.
\end{proof}
\subsection{}
\label{ranklength}
The following lemma is useful for calculations.
\begin{lem}
Let $w=s_{i_1}\ldots s_{i_t}$ be a reduced expression. Then $\ell (w)
= s(w)$ if and only if $i_j\neq i_k$ for all $j\neq k$.
\end{lem}
\begin{proof}
Recall that $s(w)$ is the minimal length of $w$ when written as a product of
arbitrary reflections. Suppose that $i_j=n=i_k$. Then
\[
s_{i_j}\ldots s_{i_k} = (s_ns_{i_{j+1}}s_n)(s_ns_{i_{j+2}}s_n)\ldots
(s_ns_{i_{k-1}}s_n) =
s_{s_n(\alpha_{i_{j+1}})}s_{s_n(\alpha_{i_{j+2}})}\ldots
s_{s_n(\alpha_{i_{k-1}})},
\]
so $\ell (w) > s(w)$.
Conversely suppose that $i_j\neq i_k$ for all $j\neq k$. We prove
that $\ell (w) = s(w)$ by induction on $\ell(w)$, the case of
$\ell(w)\leq 1$ being clear. Let $w' = s_{i_1}w$ so that
$\ell(w')=\ell (w)-1$. Suppose that $\beta=\sum n_i\alpha_i\in
Q^w$. Since a reduced expression of $w'$ does not contain the
simple reflection $s_{i_1}$ we deduce that if $w'\beta = \sum
n_i'\alpha_i$ then $n_{i_1} = n'_{i_1}$. As $\beta\in Q^w$ we have
\[
w'\beta = s_{i_1}\beta = \beta - <\beta,
\alpha^{\vee}_{i_1}>\alpha_{i_1}.
\]
As a result  $<\beta, \alpha^{\vee}_{i_1}> = 0$, implying that
$w'\beta = \beta$. Thus $P^w \subseteq P^{w'}$ and so we are in
the case where $P^w = P^{w'} \cap P^{s_i}$, (see
\cite[Section 5.3]{decpro49}), so that $s(w) = s(w') + 1 = \ell (w') +
1 = \ell(w)$. This proves the lemma.
\end{proof}
\subsection{}
We now have an upper and a lower bound on the number of blocks of
$\UP (b)$ for $b\in X_{w,e}$. Namely, if $k$ is the number of
simple modules in a block then it follows from Theorem
\ref{Upquiv} that, with $d$ as defined there, $\ell^d \leq k \leq
\ell^{r-s(w)}$. We present a sufficient condition for these bounds
to agree.
\begin{cor}
Let $b\in X_{w,e}$ where $w= s_{i_1} \ldots s_{i_t}$ is a reduced expression such that
$i_j\neq i_k$ for all $j\neq k$. Then $\UP (b)$ has a unique block.
\end{cor}
\begin{proof}
By Theorem \ref{Upquiv}(ii), under the hypothesis of the corollary
$B^w_{\ell}$ has cardinality $\ell^{r-t}$. By Lemma
\ref{ranklength}, $t= \ell (w) = s(w)$, proving the corollary, in
the light of Theorem \ref{Upquiv}(iii).
\end{proof}
\begin{rem}
Under the circumstances of the corollary, its proof together with
Theorem \ref{Upquiv}(iii) shows that the quiver of $\UP (b)$ is a
multiply-edged Cayley graph of $P_{\ell}^w$, that $P_{\ell}^w =
\sum_{\varpi_j \in \mathcal{C}_w}\bar{\mathbb{Z}}\varpi_j$, and
that the graph includes an arrow starting at each vertex for each
fundamental weight $\varpi_j$ in $\mathcal{C}_w$. But we don't
know whether additional arrows can also occur in the Cayley graph,
besides copies of these ones.
\end{rem}
\subsection{}
\label{ublocks} In many cases we can determine the number of
blocks of $\UP (b)$ for all $b \in B$ very easily.
\begin{thm}
(i) Let $b\in X_{w,e}$ and $b'\in \overline{X_{w,e}}$. Then the number of blocks of $\UP (b')$ is no greater than the number of blocks of $\UP (b)$.
\\
(ii)  For $b\in B$ the algebra $\UP (b)$ has at most $\ell^{r-s(w_0)}$
blocks, and this upper bound is attained for $b$ in the (open,
dense) stratum $X_{w_0,e}$.
\\
(iii) The following are equivalent:\\
1. $Z(U^{\geq 0}_{\epsilon}) = Z_0$;\\ 2. $s(w_0) = r;$\\ 3.
 the Cartan matrix $C$ is of type $B_r$, $C_r$, $D_r$ ($r$ even),
$E_7$, $E_8$, $F_4$ or $G_2$;\\ 4. for all $b \in B$,
 $\UP (b)$ has a unique block.
\end{thm}
\begin{proof}
(i) By \cite[Proposition 2.7]{gab} the set $\{ A\in \text{Alg}(n) : \text{number of blocks of $A \leq s$}\}$ is closed in $\text{Alg}(n)$. Thus, by Remark \ref{degen} the algebra $\UP (b')$ has no more blocks than $\UP (b)$.
\indent(ii) For $b$ in the non-empty open set $X_{w_0,e}$ of $B^-$,
$\UP (b)$ is semisimple with $\ell^{r-s(w_0)}$ simple modules, by
Proposition \ref{description}. The first part of the claim now follows from (i).
\indent(iii) $1
\Longleftrightarrow 2:$ Since $Z_0$ is integrally closed, $Z_0
\varsubsetneqq Z(\UP)$ if and only if these algebras have distinct
quotient fields, and this happens if and only if, for a generic
maximal ideal $\mathfrak{m}$ of $Z_0$, $Z/\mathfrak{m}Z
\varsupsetneqq \C$. Since the generic maximal ideal of $Z_0$ is
contained in a maximal ideal of $Z$ unramified over $Z_0$, we can
conclude that $Z_0 \varsubsetneqq Z(\UP)$ if and only if there is
a maximal ideal of $Z_0$ contained in at least two maximal ideals
of $Z$. So the equivalence now follows from Theorem
\ref{muller}.\\ $2 \Longleftrightarrow 4$: By (ii).\\ $2
\Longleftrightarrow 3$: This can be read off from \cite[Table
1]{gor4}.
\end{proof}
\subsection{Examples}
We illustrate the above analysis with a couple of examples.
(i) R = $\BARUP$. In this case the simply transitive
 group of winding automorphisms is $Q_{\ell}$. A simple calculation yields:
\begin{eqnarray*}
\tau_{\alpha_i} (E_j) = E_j \\
\tau_{\alpha_i} (K_{\lambda}) = \ep^{(\lambda , \alpha_i)}K_{\lambda}.
\end{eqnarray*}
A set of primitive idempotents of $R$ is given by
\[
e_{\mu} = \sum_{\lambda \in P_{\ell}} \ep^{(\lambda ,\mu)}
K_{\lambda},
\]
for $\mu \in Q_{\ell}$. Thus $\tau_{\alpha_i}(e_{\mu}) = e_{\mu
+\alpha_i}$. In the notation of Section 5 it is easy to
check that $R_1$ is the subalgebra of $R$ generated by the elements
$E_i$ for $1\leq i\leq r$, so is just $\BARUN$. Under the
identification of $X(Q_{\ell})$ with $P_{\ell}$ using the inner
product, it is straightforward to check that $y_{\lambda} =
K_{\lambda}$ for all $\lambda \in P_{\ell}$. Therefore the analysis of
the Section 5 recovers the well-known description of $\BARUP$
as a skew group extension of $\BARUN$. Note too that in this case $B_{\ell}^e
= P_{\ell}$ so there is a unique block. In fact one can check that
$J_1/J_1^2$ has a basis given by the images of $E_1,\ldots ,E_r$. The
quiver of $R$ (together with its relations) is described in \cite{gor2}.
(ii) $R= \UP (b)$ where $b \in X_{w,e}$ with $w=w_0s_i$. The
following descriptions can be read off from \cite[Theorem
7.7]{gor4}. We present the algebras in the format of Proposition
\ref{skewgr}, that is as matrix rings over skew group algebras
whose coefficient rings are scalar local. There are two cases,
depending on the value of $w_0(\alpha_i)$. In the ring-theoretic
context of Proposition \ref{blquiv} the dichotomy is determined by
whether or not $\tilde{P}_{\ell}^w$ acts trivially on the Jacobson
radical.
\\
(a) $w_0 (\alpha_i) = -\alpha_i$. Here,
 \begin{eqnarray*}
\UP (b) \quad \cong \quad
\mathrm{Mat}_{\ell^{1/2(N-1+s(w))}}(\BARUP(\mathfrak{sl}_2) \otimes
\C C_{\ell}^w),
\end{eqnarray*}
where $|C_{\ell}^w| = {\ell}^{r-s(w)-1} = {\ell}^{r-s(w_0)}$ and we note that $\BARUP(\mathfrak{sl}_2)$ is a skew group extension of a basic algebra as in (i).\\
 (b) $w_0 (\alpha_i ) \neq -\alpha_i$. In this case,
 \begin{eqnarray*}
\UP (b) \quad \cong \quad
\mathrm{Mat}_{\ell^{1/2(N-1+s(w))}}(\C[X]/<X^{\ell}> \otimes \C
C_{\ell}^w),
 \end{eqnarray*}
 where $C_{\ell}^w = \tilde{P}_{\ell}^w$, and $|\tilde{P}_{\ell}^w| = {\ell}^{r-s(w)} = {\ell}^{r-s(w_0)-1}$.
\section{Reduced quantised function algebras}
\subsection{}
Recall the definition of $\bi, \ci \in \EO$ given in \ref{centre}.
\label{firstcount}
\begin{lem}
Let $g\in G$. The number of blocks
in $\EO (g)$ equals $\ell^d$ where $d$ is the cardinality of the set
$\{ 1\leq i\leq r: \bi^{\ell}(g) \neq 0 \neq \ci^{\ell}(g)\}$.
\end{lem}
\begin{proof}
We freely use the notation of the proof of Theorem \ref{azmor}. By
Theorem \ref{muller} it is sufficient to show that $Z_g$ has exactly
$d$ maximal ideals. By the proof of Theorem \ref{azmor}
\[
Z_g \cong \bigotimes_{i=1}^r R_i,
\]
where $R_i$ is as in (\ref{jess}). Let $b_i = \bi^{\ell}(g)$ and
$c_i = \ci^{\ell}(g)$. We have already seen in (\ref{eq:ss}) that
if $b_i\neq 0 \neq c_i$ then $R_i$ is isomorphic to
$\mathbb{C}^{\ell}$, whilst if $b_i\neq 0$ and $c_i = 0$ then
$R_i$ is isomorphic to $\mathbb{C}[X]/(X^{\ell})$ by
(\ref{eq:trun}). It remains to consider the case $b_i = 0 = c_i$.
In this case inclusion provides an isomorphism
\[
\frac{\C [X_1, \ldots ,X_{\ell -1}]}{I_i'} \longrightarrow R_i
\]
where $I_i'$ is the ideal generated by $X_kX_{k'}$ for all $k$ and
$k'$. Since this is manifestly a local ring, the result follows
from these calculations and Theorem \ref{muller}.
\end{proof}
\subsection{}
As noted in \cite[Appendix]{declyu1} the elements $\bi^{\ell},
\ci^{\ell}\in \mathcal{O}[G]$ can be interpreted as the classical
analogues of the quantum matrix coefficients $\bi$ and $\ci$. In
other words $\bi^{\ell}$ (respectively $\ci^{\ell}$) is the
classical matrix coefficient
$\tilde{c}^{\varpi_i}_{f_{-w_0\varpi_i},v_{\varpi_i}}$
(respectively
$\tilde{c}^{-w_0\varpi_i}_{f_{w_0\varpi_i},v_{-\varpi_i}}$). Here
we have used $\tilde{c}$ to distinguish the classical matrix
coefficients from their quantum analogues. \label{Weylcomb}
\begin{lem}
Let $g\in X_{w_1,w_2}$. The algebra $\EO (g)$ has exactly $\ell
^d$ blocks where $d$ is the cardinality of the set $\{ 1\leq i\leq
r: w_0w_1,w_0w_2 \in \mathrm{Stab}_W(\varpi_i)\}$.
\end{lem}
\begin{proof}
Let's write $B_i \in \mathcal{O}[G]$ for the regular function
$\bi^{\ell}$, and similarly $C_i$ for $\ci^{\ell}$. Let $w$ be an
arbitrary element of $W$. We have
$f_{-w_0\varpi_i}(bwb'v_{\varpi_i}) \in
\mathbb{C}^*f_{-w_0\varpi_i}(bwv_{\varpi_i})$ for $b,b'\in B$ since
$v_{\varpi_i}$ is a highest weight vector. Hence $B_i(bwb') =
f_{-w_0\varpi_i}(bwb'v_{\varpi_i})$ is not identically zero
if and only if $wv_{\varpi_i}$ is a lowest weight vector, and in
this case it is non-zero for all $b,b' \in B$. Note that this
happens if and only if  $w\varpi_i = w_0\varpi_i$. A similar
analysis applied to $C_i$ shows $C_i (bwb') \neq 0$ for some $b,b'
\in B$ if and only if $C_i (bwb') \neq 0$ for all $b,b' \in B$ if
and only if $wv_{-\varpi_i}$ is a highest weight vector. That is,
 if and
only if $-w\varpi_i=-w_0\varpi_i$.
We conclude from the previous paragraph that if $g\in X_{w_1,w_2}$
then $B_i(g)\neq 0$ and $C_i(g)\neq 0$ if and only if
$w_0w_1,w_0w_2\in \mathrm{Stab}_W(\varpi_i)$. The lemma now
follows from Lemma \ref{firstcount}.
\end{proof}
\subsection{}
\label{finalversion} Recall the definition of the winding
automorphisms given in \ref{appnot}. For the rest of the paper we
will fix $w_1, w_2 \in W$ and assume that $\ell$ is prime to
$\text{ord}(w_2^{-1}w_1)$.
\begin{thm}
Let $g\in X_{w_1,w_2}$ and let $w=w_2^{-1}w_1$. Let
$N(w_1,w_2)\subseteq Q_{\ell}^w$ be the normaliser of one (hence
any) block of $\EO (g)$ with respect to the (right) winding
automorphisms. Let $\mathfrak{S}(w_1,w_2) = \{1\leq i\leq
r:w_0w_1,w_0w_2\in \mathrm{Stab}_W(\varpi_i)\}$. Then
\[
N(w_1,w_2) = Q_{\ell}^w \cap (\sum_{i\notin
\mathfrak{S}(w_1,w_2)}\overline{ \mathbb{Z}\alpha_i}).
\]
\end{thm}
\begin{proof}
By Theorem \ref{muller} different blocks arise from the different
maximal ideals of $Z$ lying over $\mathfrak{m}_g$, so we need to
see how these are permuted by the right winding automorphisms. By
\ref{centre} such maximal ideals are determined by the central
elements $\bi^k\ci^{\ell -k}$ for $1\leq i\leq r, \, 0 \leq k
\leq \ell$. Hence we need only study the action of the winding
automorphisms on $\bi$ and $\ci$.
By definition $\bi = c^{\varpi_i}_{f_{-w_0\varpi_i},v_{\varpi_i}}$ where both $f_{-w_0\varpi_i}$ and $v_{\varpi_i}$ are
highest weight vectors. Therefore
\[
\Delta (\bi) = \sum_{j} c^{\varpi_i}_{f_{-w_0\varpi_i}, v_j}\otimes
c^{\varpi_i}_{f_j, v_{\varpi_i}}
\]
where $\{f_j\}$ and $\{v_j\}$ are dual bases of $V(\varpi_i)^*$ and
$V(\varpi_i)$ respectively. Let $\beta \in Q_{\ell}$. We have
\[
\beta(c^{\varpi_i}_{f_{-w_0\varpi_i},v_j}) =
\begin{cases}
\ep^{(\beta , \varpi_i)} \qquad &\text{if $f_{-w_0\varpi_i}$ and $v_j$ are dual,}\\
0& \text{otherwise.}
\end{cases}
\]
Therefore, letting $\tau_{\beta}$ denote the right winding
automorphism of $\EO$ defined analogously to those for $\UP$ in
\ref{appnot}, $\tau_{\beta}(\bi) = \ep^{(\beta, \varpi_i)}\bi$.
Similarly, we find that $\tau_{\beta}(\ci) = \ep^{-(\beta ,
\varpi_i)}\ci$. We have shown that
\[
\tau_{\beta}(\bi^k\ci^{\ell -k}) = \ep^{2k(\beta ,
\varpi_i)}\bi^k\ci^{\ell - k}.
\]
As in the proof of Theorem \ref{azmor} the maximal ideals of $Z$ lying over $\mathfrak{m}_g$ are obtained by
piecing together the maximal ideals of the algebras $R_i$ for $1\leq
i\leq r$ (notation of the proof of Theorem \ref{azmor}). These maximal ideals in turn depend on
the vanishing behaviour of $\bi^{\ell}(g)$ and $\ci^{\ell}(g)$. Specifically if
$\bi^{\ell}(g)=0$ or $\ci^{\ell}(g) = 0$ then $R_i$ has a unique
maximal ideal whilst if $\bi^{\ell}(g)\neq 0 \neq \ci^{\ell}(g)$ then
$R_i$ is semisimple with exactly $\ell$ maximal ideals. In the
second case it follows from the previous paragraph that the winding
automorphisms permute the primitive idempotents of $R_i$
non-trivially unless $(\beta ,\varpi_i) \in \ell\mathbb{Z}$. Hence
the same is true of the maximal ideals. We deduce that the
normaliser in $Q_{\ell}$ of a block of $\EO (g)$ is simply the
subgroup
\[
\{ \beta \in Q_{\ell} : (\beta , \varpi_i) \in \ell\mathbb{Z}
\textit{ whenever } \bi^{\ell}(g) \neq 0 \neq \ci^{\ell}(g) \}.
\]
By the proof of Lemma \ref{Weylcomb} this equals the subgroup
\[
\{ \beta \in Q_{\ell} : (\beta , \varpi_i ) \in \ell \mathbb{Z}
\textit{ whenever } w_0w_1, w_0w_2 \in \mathrm{Stab}_W(\varpi _i)\}.
\]
The theorem follows from Lemma \ref{ellform}(iii).
\end{proof}
\subsection{}\label{oblocks}
Let $g\in X_{w_1,w_2}$ and let $w=w_2^{-1}w_1$. The factor group
$Q^w_{\ell}/N(w_1,w_2)$ acts simply transitively on the blocks of
$\EO (g)$. We claim that this factor group is an elementary
abelian $\ell$-group of rank the cardinality of
$\mathfrak{S}(w_1,w_2)$, where $\mathfrak{S}(w_1,w_2)$ is as in
the statement of Theorem \ref{finalversion}. To see this recall
from the proof of Lemma \ref{ellform} that we have a
$<w>$-invariant decomposition $P_{\ell} = P^w_{\ell} \oplus
P'_{\ell}$. Moreover, since $\ell$ is prime to the order of $w$,
we have $Q^w_{\ell}= (P'_{\ell})^{\perp}$. By definition there is
an isomorphism
\[
\frac{Q^w_{\ell}}{N(w_1,w_2)} \cong \frac{Q^w_{\ell} + \sum_{i\notin
\mathfrak{S}(w_1,w_2)}\overline{\mathbb{Z}\alpha_i}}{ \sum_{i\notin
\mathfrak{S}(w_1,w_2)}\overline{\mathbb{Z}\alpha_i}}.
\]
If $i\in \mathfrak{S}(w_1,w_2)$ then $w\varpi_i =
(w_0w_2)^{-1}(w_0w_1)\varpi_i =\varpi_i$. Thus $ \sum_{i\in
\mathfrak{S}(w_1,w_2)}\overline{\mathbb{Z}\varpi_i} \subseteq
P^w_{\ell}$ which in turn implies that $ \sum_{i\in
\mathfrak{S}(w_1,w_2)}\overline{\mathbb{Z}\varpi_i}\cap P'_{\ell} = 0$. We deduce
that
\[
(Q^w_{\ell} +  \sum_{i\notin
\mathfrak{S}(w_1,w_2)}\overline{\mathbb{Z}\alpha_i})^{\perp} =
(Q^w_{\ell})^{\perp} \cap ( \sum_{i\notin
\mathfrak{S}(w_1,w_2)}\overline{\mathbb{Z}\alpha_i})^{\perp} = P'_{\ell}
\cap  \sum_{i\in
\mathfrak{S}(w_1,w_2)}\overline{\mathbb{Z}\varpi_i} = 0.
\]
Therefore $Q^w_{\ell} +  \sum_{i\notin
\mathfrak{S}(w_1,w_2)}\overline{\mathbb{Z}\alpha_i} = Q_{\ell}$ and we
see that the factor group is isomorphic to
\[
\frac{Q_{\ell}}{ \sum_{i\notin
\mathfrak{S}(w_1,w_2)}\overline{\mathbb{Z}\alpha_i}},
\]
an elementary abelian $\ell$-group of the required rank. Summing
up, we have shown:
\begin{cor}
 Let $g \in X_{w_1,w_2}$. Write $w = w_2^{-1}w_1$.\\
 1. $\EO (g)$
has $\ell^{\mathrm{card}(\mathfrak{S}(w_1,w_2))}$ blocks, each
containing $\ell^{r - s(w) -
\mathrm{card}(\mathfrak{S}(w_1,w_2))}$ simple modules.\\ 2. The
quiver of each block of $\EO (g)$ is a multiply-edged Cayley graph
of $N(w_1,w_2)$.
\end{cor}
\subsection{Examples}
\label{oex}
 (i) Suppose first that $g\in X_{e,e}$. Here we have
$\mathfrak{S}(e,e) = \emptyset$ and so $N(e,e) = Q_{\ell}$. In
particular there is only one block. The quiver (and the relations)
are described in \cite{gor2}.
\\
(ii) Suppose that $g\in G$ lies on the Azumaya locus. Recall from
(\ref{maxfactor}) that this can be expressed by the condition that
$g\in X_{w_1,w_2}$ where $\ell(w_1) + \ell(w_2) +s(w_2^{-1}w_1)
=2N$. Let $u = w_0w_1$, $v=w_0w_2$ and $w = v^{-1}u =
w_2^{-1}w_1$. Thus $\ell (u) + \ell(v) = s(w)$. We have
inequalities
\[
\ell(w) \leq \ell(v^{-1}) + \ell(u) = s(w) \leq \ell(w),
\]
which must be equalities. In particular $\ell (w) =
s(w)$ so by Lemma \ref{ranklength} $w= s_{i_1}\ldots s_{i_t}$
where $i_j\neq i_k $ for all $j \neq k$. Since $\ell(v) + \ell(u)
= \ell(w)$ we can assume without loss of generality that $v =
s_{i_k}\ldots s_{i_1}$ and $u=s_{i_{k+1}}\ldots s_{i_t}$. It is
now clear that
\[
\{ 1\leq i \leq r : w_0w_1 \in \mathrm{Stab}_W(\varpi_i)\} =
I\setminus \{i_{k+1},\ldots ,i_t\},
\]
and that
\[
\{ 1\leq i \leq r : w_0w_2 \in \mathrm{Stab}_W(\varpi_i)\} =
I\setminus \{i_1,\ldots ,i_k\}.
\]
Hence $\mathfrak{S}(w_1,w_2) = I\setminus\{i_1,\ldots ,i_t\}$ has
cardinality $r-t = r - s(w)$. Thus we have a group of order
$\ell^{r-s(w)}$ permuting the blocks simply transitively and a unique
simple in each block. This agrees with Theorem \ref{azmor}.
\\
(iii) Let $w_1=w_0s_i=w_2$. In this case $\mathfrak{S}(w_1,w_2) = I\setminus
\{ i\}$. Thus $\EO(g)$ has $\ell^{r-1}$ blocks and each contains
$\ell^{r-(r-1)} = \ell$ simple modules.
\\
(iv) Let $G= SL_4(\mathbb{C})$. We will compare two
situations:
\\
\indent
(a) $w_0w_1=s_1s_2s_3$, $w_0w_2=s_1$: here $\mathfrak{S}(w_1,w_2)= \emptyset$ ;
\\
\indent
(b) $w_0w_1= s_1s_2s_1$, $w_0w_2=s_1$: here $\mathfrak{S}(w_1,w_2) = \{ 3 \}$ .
\\
In both cases we have $\ell (w_1) = 7$, $\ell (w_2) =9$ and $s(w_2^{-1}w_1) =
2$. This shows that $\ell(w_1), \ell(w_2)$ and $s(w_2^{-1}w_1)$ are not a
complete set of invariants for the algebras $\EO (g)$. This was
unknown before!

\end{document}